\documentclass[11pt]{article}
\pdfoutput=1
\usepackage{algorithmic,algorithm}
\usepackage{float}
\usepackage{enumitem}
\usepackage{multirow}
\usepackage{latexsym,amssymb,amsmath,amsthm,graphicx,float}
\usepackage{bm}
\usepackage{natbib}

\newcommand{\beq}{\begin{equation}}
\newcommand{\eeq}{\end{equation}}
\renewcommand{\tilde}{\widetilde}
\renewcommand{\hat}{\widehat}
\newcommand{\bit}{\begin{itemize}}
\newcommand{\eit}{\end{itemize}}
\newcommand{\ben}{\begin{enumerate}}
\newcommand{\een}{\end{enumerate}}

\title{Velocity estimation via registration-guided \\ least-squares inversion}
\author{ Hyoungsu Baek$^1$, Henri Calandra$^2$, Laurent Demanet$^1$\\
$^1$Department of Mathematics, MIT, Cambridge, MA 02139, USA \\
$^2$TOTAL Exploration \& Production, Avenue Larribau, 64018 Pau, France }
\date{March 2013}		% remove for today's date

\begin{document}

\maketitle

\vspace{-1cm}
\begin{abstract}
This paper introduces an iterative scheme for acoustic model inversion where the notion of proximity of two traces is not the usual
least-squares distance, but instead involves registration as in image processing. 
Observed data are matched to predicted waveforms via piecewise-polynomial warpings, obtained by solving a nonconvex optimization problem in
a multiscale fashion from low to high frequencies. 
This multiscale process requires defining low-frequency augmented signals in order to seed 
the frequency sweep at zero frequency.  Custom adjoint sources are then defined from the warped waveforms. 
The proposed velocity updates are obtained as the migration of these adjoint sources, and cannot be interpreted as the negative gradient of any given objective function.
The new method, referred to as RGLS, is successfully applied to a few scenarios of model velocity
estimation in the transmission setting.  We show that the new method can converge to the correct model in situations where conventional
least-squares inversion suffers from cycle-skipping and converges to a spurious model.
\end{abstract}

{\bf Acknowledgments.} 
The authors are grateful to the authors of Madagascar; and to TOTAL S.A. for support. 

\section{Introduction}

Waveform inversion via non-linear least-squares (LS) minimization
\cite[]{Tarantola1982} is effective when the starting model is accurate \cite[]{Virieux2009}, 
but otherwise suffers from stalled convergence to spurious local minima due to the oscillatory nature of the data. Traveltime tomography is
typically used to generate an accurate initial model \cite[]{GPR1099,Pratt01021991,Bregman01021989}, which is then improved upon by waveform inversion.  
Frequency sweeps in full waveform inversion (FWI) \cite[]{Pratt99seismicwaveform} are another 
well-documented way to encourage convergence toward the global least-squares
minimum, by fitting data from low to high frequencies.  However, the lack of low-frequency data, or their corruption by noise, often hinders
this frequency sweep approach. This paper proposes a method which succeeds at recovering the model velocity in some transmission scenarios,
without traveltime picking, and even when the data only contain high frequencies. 

We propose a simple modification to least-squares-based waveform inversion where the residual of the  adjoint-state equation, normally the
difference $d-u$ between the observed data $d$ and predicted data $u$, is replaced by a more geometrically meaningful quantity.
We propose
to let this more general residual go by the name {\em adjoint source}, and form it by replacing the observed data $d$ with a transported, or
warped, version $\widetilde{d}$ of the predicted data toward the observed data. The rationale for this substitution is that the phases of
$d$ are in general off by more than one wavelength in comparison to those of $u$. In contrast, the new {\em fractionally warped data}
$\widetilde{d}$ are defined so that their phases match those of the prediction $u$ to within a small fraction of a wavelength. When the
registration algorithm generating the warping succeeds, this method yields model velocity updates that properly fix time delays or advances without suffering from cycle skipping --- although these velocity updates are no longer gradient descent steps.
 We observe that least-squares misfit optimization with this modified adjoint source can have a much
enlarged basin of attraction. The method is referred to as {\em registration-guided least-squares} (RGLS) inversion.

To find the best warping between an observed trace $d$ and the corresponding predicted trace $u$,
we formulate and solve an optimization problem which, not unlike FWI, is itself nonconvex. The highly oscillatory nature of the traces is
also what makes seismogram registration nontrivial. 
However, we show that the nonconvexity is tractable and can be handled by a continuation strategy, where the match is realized
scale-by-scale in a careful, iterative fashion. The traces $d$ and $u$ usually do not contain useful low frequencies in exploration
seismology, so the seeding problem of this multiscale iteration is as much an issue here as in classical FWI. We propose to solve this
problem by introducing nonlinearly transformed signals which, by construction, contain low-frequency components. We refer to these
convenient, nonphysical nonlinearly-transformed signals as {\em low-frequency augmented} (LFA) signals. The LFA transformation can be thought
of as an ad-hoc pre-processing of the traces so as to create low frequencies, yet maintain much of the information at high frequencies.
Subsequently, seismogram registration is realized through the match of the LFA of $d$ to the LFA of $u$ by a warping function of limited
complexity, such as a piecewise polynomial.

Matching observed data and predicted data, or matching time-lapse images 
has already been used in many different applications under many different names: 
phase-shift \cite[]{Tromp2005}, traveltime delay based on correlation \cite[]{luo645}
image registration using optimal transport \cite[]{938878}, curve registration \cite[]{ramsay1998}, 
registration using local similarities \cite[]{FomelJin2009,FomelBaan2010}, 
and dynamic time warping for speech pattern matching \cite[]{saoke78}.
\cite{Anderson1983221} applied the dynamic waveform matching method 
to geophysics problems such as earthquake arrival time identification and waveform clustering.
\cite{HaleGeoPhysics2013} applied the dynamic warping, which was originally developed for speech recognition, to analysis of time-lapsed
seismic images.
\cite{Maggi2009} developed an algorithm to select time windows
to extract time-shift between seismograms with iterative tomographic inversion in mind.
\cite{Liner01112004} studied alignment of seismic traces using dynamic programming 
for trace processing and interpretation.
\cite{Kennett2012GJI} introduced a transfer operator which maps seismograms
 similar to our polynomial mapping as a way to generalize comparison between traces.
 
It is important to point out that the ``fractional" character of the warping used to generate the adjoint source in this paper is in the
same spirit as a solution proposed in \cite{sava_biondi2004geoprospect}, where image perturbations in the image domain are replaced by their
linearized version to mitigate lack-of-convexity issues beyond the Born approximation.  
Similar suggestions of non-gradient updates generated from residuals involving `infinitesimally modified" images are also proposed in recent work by \cite{Fei2010} and \cite{Albertin2011}.

To the best of our knowledge, no studies have yet proposed to modify the adjoint source by replacing 
observed data with time-warped predicted data in order to enlarge the basin of attraction of FWI.
Moreover, we hope to illustrate the benefits of considering piecewise polynomials to define mappings between two images or traces. Finally,
the idea of nonlinearly transforming traces to generate low frequencies for seismogram registration seems new to seismic imaging.

The paper is organized as follows. We start by explaining the motivation behind 
modifying the adjoint source to the adjoint state equation. We then detail  
and compare different trace registration methods.
Seismogram registration at the trace level is demonstrated with 
synthetic noisy and noiseless data. The RGLS inversion method is then tested in several  
transmission cases, with velocity models that include low and high velocity zones. 
We show that the LS and RGLS methods behave significantly differently.
We finish by discussing the limits of the proposed method.
In a nutshell, seismogram registration requires comparable traces, which explains why we consider transmission rather than reflection
examples in this paper.

\section{Guided least-squares with a modified adjoint source}
\label{sec:modifiedLSQ}

Full waveform inversion (FWI), in its standard form, tries to minimize the least-squares misfit 
\begin{equation}
J[m]=\frac{1}{2} \sum_{s,r} \int | u_s (x_r,t) - d_s (x_r,t)|^2 dt,
\label{eq:least_square_misfit}
\end{equation}
 where $m(x)$ denotes a velocity model,
$u_s= {\cal F}_s[m]$, and $d_s$ are predicted and observed data at a shot $s$ and at a receiver $x_r$,
respectively. For notational simplicity, the subscripts $s$ and $r$ are omitted whenever it
does not cause confusion. In this paper the forward operator ${\cal F}_s[m]$ maps a velocity model $m$ to data $u_s(x_r,t)$ through the
acoustic wave equation, 
\[
m \frac{\partial^2 u_s}{\partial t^2} = \Delta u_s+f_s (x,t),
\] 
where $f_s(x,t)$ is a source  term. The adjoint-state method generates the gradient of $J[m]$ via the
imaging condition
\begin{equation}
\frac{\delta J}{\delta m}[m] = \int q_s (x,t)  \frac{\partial^2  u_s}{\partial t^2} (x,t) dt,
\end{equation}
where the adjoint field $q_s$ solves the adjoint wave equation backward in time with the data residual
$d_s - u_s$  in the right-hand side, and in the medium $m$.

\begin{figure}
\begin{center}
\hspace{-0.5cm}
\includegraphics[width=1\textwidth]{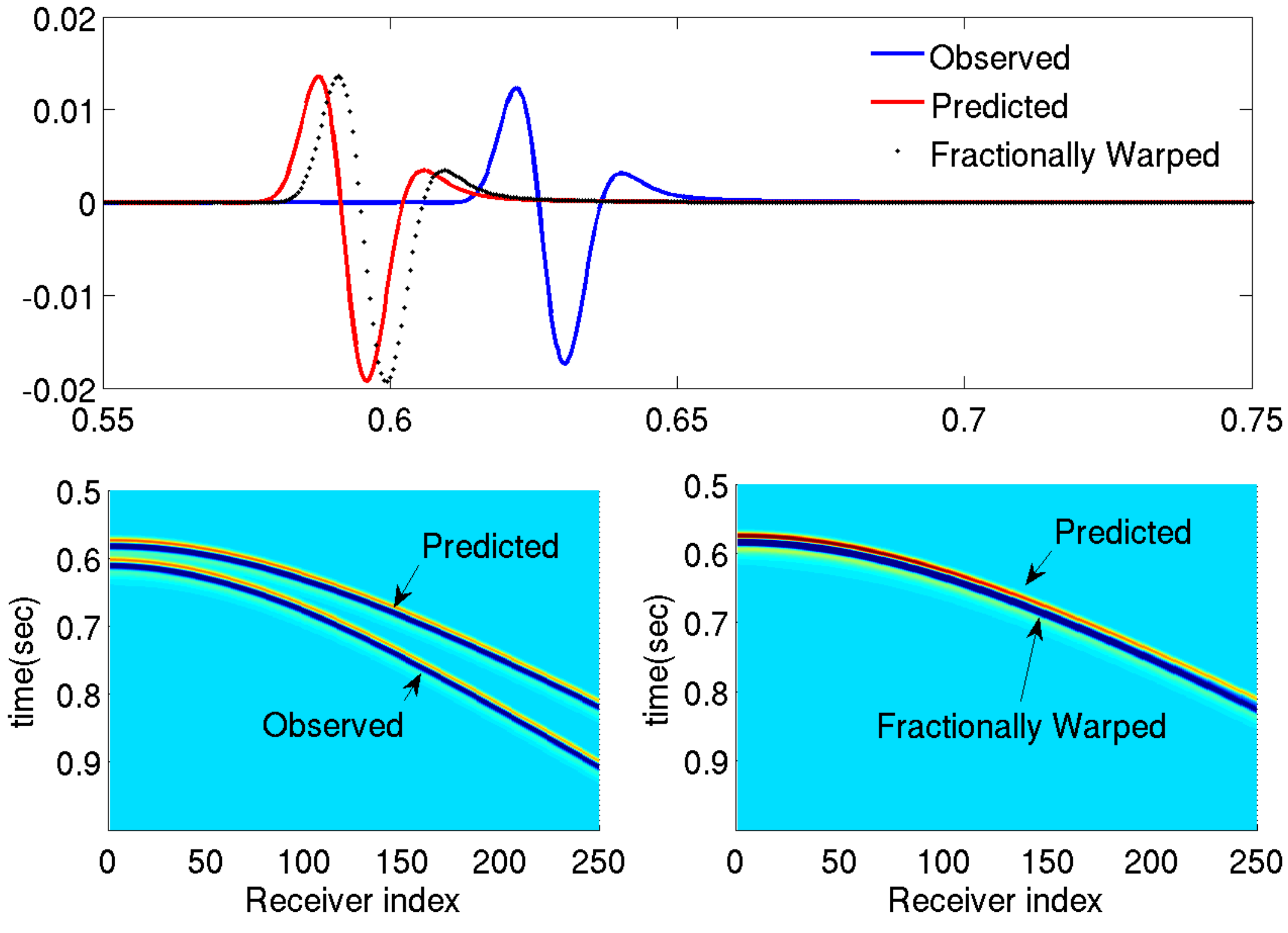} %%Fractionally_warped.pdf}
\end{center}
\caption{
Replacement of observed data with fractionally warped data, which are a  mapped (transported slightly) version of the given predicted data
towards the observed data. We call the new traces {\it fractionally warped data}.
}
\label{fig:fictitious_observed_data}
\end{figure}

It is well known that the nonconvexity of $J[m]$ is particularly pronounced way when the data are oscillatory. More specifically, when the
time difference between corresponding arrivals in $u_s$ and $d_s$ is
larger than a half wavelength, the steepest descent direction of the data misfit may
result in increasing those time differences, consequently increasing the model error. In
order to guide the optimization in the direction of the correct model update, we propose
to change the residual $d_s-u_s$ in the adjoint wave equation by replacing $d_s$ by a version
of $u_s$ transported ``part of the way" toward $d_s$. We denote these virtual, transported data by $\tilde{d}_s$ and
refer to them as {\em fractionally warped data}. Their construction is in the next section. The
rationale behind $\tilde{d}_s$ is that its arrivals can now be less than a quarter of a wavelength
apart from those in $u_s$. 
The substitution of $d_s$ by $\tilde{d}_s$ is illustrated in Figure \ref{fig:fictitious_observed_data}.

In order to generate a good candidate of fractionally warped data
we propose to find piecewise cubic polynomials $p(t)$ and $A(t)$ 
so as to have a good match $d(t) \approx A(t) u(p(t))$, then define fractionally warped data as
\begin{equation}
\tilde{d}(t) = \left[ A(t) \right]^{\alpha} u((1-\alpha)t + \alpha p(t))
\end{equation}
with some very small $0<\alpha \ll 1$.
The warping $p(t)$, the amplitude $A(t)$, and the value of $\alpha$ are chosen
so that fractionally warped data $\tilde{d}(t)$ have a similar shape to that of the prediction $u(t)$ but
phase discrepancies smaller than a quarter of a wave period.
In order to find $p(t)$ and $A(t)$, we propose a non-convex optimization scheme similar to image registration.
The proposed FWI method therefore transfers all non-convexity  
to the registration problem at the trace level. 
We describe seismogram registration in detail in the next section.

\section{Seismogram registration}
\label{sec:registration}

\subsection{Statement of the optimization problem}
In order to find $A(t)$ and $p(t)$ we propose to solve the following least-squares
minimization problem for each trace: find $p(t)$ and $A(t)$ piecewise cubics that minimize
\begin{equation}
W[p,A] = \frac{1}{2} \int |d(t) - A(t) u(p(t)) |^2 dt + \frac{1}{2} \int |p(t)-t|^2 dt. 
\label{eq:opt1}
\end{equation}
The letter $W$ stands for ``warping". The functions $A(t)$ and $p(t)$ are realized as cubic Hermite splines over 4 fixed subintervals in the
following examples, unless otherwise stated. (This paper does not address the problem where the spline nodes could also be determined by
optimization.)

Due to the oscillatory nature of the predicted and observed data, 
this registration problem is non-convex and suffers from the same 
cycle-skipping phenomenon as conventional least-squares FWI does.
Simulated annealing or other Monte Carlo methods for global optimization 
could be performed but they are not tried in this paper.
Instead, the minimization is carried out in a multiscale fashion 
by restricting the data $d(t)$ and its prediction $A(t) u(p(t))$ to a slowly growing subset of frequencies, from the DC component to
successively higher frequencies, as in frequency domain FWI \cite[]{Plessix2009}. 
However, the observed data $d$ usually have small energy in the low frequency band and may be corrupted by noise. In order to start the
sweep at zero frequency, 
we use modified traces $D(t)$ and $U(t)$, manufactured to contain low frequencies, instead of $d(t)$ and $u(t)$. We call $D(t)$ and $U(t)$
{\em low-frequency augmented} (LFA) signals. Hence, the optimization problem (\ref{eq:opt1}) becomes: find $p(t)$ and $A(t)$ piecewise
cubics that minimize
\begin{equation}
W_{\mbox{\scriptsize LFA}} [p,A] =  \frac{1}{2} \int |D(t) - A(t) U(p(t)) |^2 dt + \frac{1}{2} \int |p(t)-t|^2 dt. 
\label{eq:opt2}
\end{equation}

Here are three reasonable possibilities for defining an LFA signal, $U(t)$, from a band-limited signal $u(t)$:
\begin{subequations}
\begin{align}
U_h &= u(t) + | u(t) + i ({\cal H} u) (t)|,  \label{eq:uh} \\
U_s &= u^2(t),   \label{eq:us} \\
U_a &= |u(t)|.   \label{eq:ua}
\end{align}
\end{subequations}
We use $I$ for the identity map and ${\cal H}$ for the Hilbert transform, defined in the frequency domain as
$\hat{{\cal H} u}(\omega) = - i \, \mbox{sgn}(\omega) \hat{u} (\omega)$ 
\cite[]{HarmonicAnalysis1997}, where $\mbox{sgn}(\cdot)$ is the sign function.

The Hilbert transform completes any real signal with an imaginary part, so that $u + i {\cal H} u$ is an ``analytic" signal in the sense of
having no negative frequency component.
The amplitude $\sqrt{ u^2(t) + ({\cal H} u(t))^2 }$ has the interpretation of an envelope for the signal $u + i {\cal H} u$. The Hilbert
transform is a classical tool in signal processing; it is typically used for demodulation in seismic inversion
\cite[]{Misfit_phase_envelope}.  In a later section we argue that $U_h$ is a particularly good LFA transformation.

The piecewise cubic polynomials $p(t)$ and $A(t)$ can be written as
\[
 p(t) = \displaystyle \sum_{k=0}^{n} \rho_{k} \phi_{k}(t) 
 = [\phi_{0}(t) \,\, \phi_{1}(t) ... \phi_{n}(t)][\rho_0 \,\, \rho_1 \,\, ... \,\, \rho_n ]^T  
\]
 and 
 \[
 A(t) = \displaystyle \sum_{k=0}^{n} \alpha_{k} \phi_{k}(t) 
 = [\phi_{0}(t) \,\, \phi_{1}(t) ... \phi_{n}(t)][\alpha_0 \,\, \alpha_1 \,\, ... \,\, \alpha_n ]^T,
 \]
 respectively, 
 for a set of basis functions $\phi_k (t), k=0,1,2,...,n$.
The column vectors $[\rho_0, \rho_1, ..., \rho_{n}]^T$ and $[\alpha_0, \alpha_1, ..., \alpha_{n}]^T$  
are denoted by ${\bm \rho}$ and ${\bm \alpha}$, respectively. Their components are the $2(n+1)$ parameters to be determined per trace. 

The gradient and the Hessian matrix of $W_{\mbox{\scriptsize LFA}}$ with respect to ${\bm \rho}$ and ${\bm \alpha}$ can be found
analytically, e.g.,
\begin{eqnarray}
\frac{\partial W_{\mbox{\scriptsize LFA}}}{\partial \rho_i} &=& \int  [D - A U(p)][-AU'(p) \phi_i]   +  (p - t) \phi_i  dt
\label{eq:gradient}, \\
(\mathbf{H}_{\bm{\rho}})_{ij} = \frac{\partial^2 W_{\mbox{\scriptsize LFA}}}{\partial \rho_i \partial \rho_j } &=& \displaystyle \int  
\left[ [AU'(p)]^2   
                                   - [D - A U(p)][A U''(p)] 
                                   + 1 \right] \phi_i \phi_j  dt  \label{eq:Hessian}.
\end{eqnarray}
Similar expressions can be derived for 
the gradient and Hessian with respect to the coefficients of $A(t)$.
The integrals are computed approximately using the trapezoidal rule as a quadrature.

\subsection{Optimization strategy and algorithm for seismogram registration}
We detail a way to resolve the non-convexity issue
of the optimization problem (\ref{eq:opt1}) or 
(\ref{eq:opt2}) for seismogram registration in {\bf Algorithm 1}.

The collection of discrete frequencies from 0 to $\omega_i$ is denoted by $\Omega_i = [0, \omega_i]$.
We create $M$ such sets, $\Omega_1, \Omega_2, ... , \Omega_M$, where
$\Omega_1 \subset \Omega_2 \subset ... \subset \Omega_M$ and $\omega_1 < \omega_2 < ... < \omega_M$. Below, $LPF_k ( \cdot )$ denotes the
application 
of a low-pass filter with passband $\Omega_k = \left[ 0, \omega_k \right]$. At the $k^{th}$ outer iteration step, both LFA traces $D$ and
$U$ are low-pass filtered to the frequency range $\omega \in \Omega_k$, resulting in the LF signals $D_k$ and $U_k$. We let
$W_{\mbox{\scriptsize LFA,k}}$ for the expression of $W_{\mbox{\scriptsize LFA}}$ with $D_k$ and $U_k$ in place of $D$ and $U$.

\begin{algorithm}[h]\label{algo1}
\caption{Seismogram registration. }
\begin{algorithmic}
\medskip
\STATE \textbf{input:} traces $u(t)$ and $d(t)$
\STATE \textbf{initialize:} $p(t) = t$, $A(t) = 1$
\STATE \textbf{LFA:}  $D(t)  \leftarrow LFA( d(t) ), \; U(t)  \leftarrow LFA ( u(t) )$ 

\FOR {$k = 1, 2, ..., M$ }
\STATE
\begin{tabular}{ll}
\textbf{filter:} &  $D_k(t)  \leftarrow LPF_k ( D(t) ),  \;U_k(t)  \leftarrow LPF_k ( U(t) ) $ 
\end{tabular}

\medskip
\WHILE{not converged}
\STATE
\begin{tabular}{ll}
\textbf{Compute:} & $\frac{\partial W_{\mbox{\scriptsize LFA,k}}}{\partial \bm{\rho}}$, $\frac{\partial W_{\mbox{\scriptsize
LFA,k}}}{\partial \bm{\alpha}}$, and the Hessians $\mathbf{H}_{\bm{\rho}}$,  $\mathbf{H}_{\bm{\alpha}}$ of the functional
$W_{\mbox{\scriptsize LFA,k}}$.
 \\
\textbf{Newton step:} & $\bm{\rho} \leftarrow \bm{\rho} - \mathbf{H}_{\bm{\rho}}^{-1} \frac{\partial W_{\mbox{\scriptsize LFA,k}}}{\partial
\bm{\rho}} $, $\;$
                        $\bm{\alpha} \leftarrow \bm{\alpha} - \mathbf{H}_{\bm{\alpha}}^{-1} \frac{\partial W_{\mbox{\scriptsize
LFA,k}}}{\partial \bm{\alpha}} $ 
\end{tabular}
\ENDWHILE
\ENDFOR
\STATE \textbf{output:} $p(t)$, $A(t)$
\end{algorithmic}
\end{algorithm}

The maximum frequency $\omega_M$ in the outer loop is set to a frequency below the central frequency
of the source signature used to generate data.
Numerical experiments show that sweeping up to half of the central frequency suffices for convergence. 
Most kinematic discrepancies between the observed data and the predicted data 
are fixed after sweeping up to 5 Hz, 
where central frequencies of the Ricker wavelets are above 15 Hz up to 50 Hz.
Since registering every trace is time-consuming, registration is only performed once every 50 traces. 
The mappings for the other traces are interpolated.

\subsection{Examples of registration of synthetic traces} 

\begin{figure}
\begin{center}
\includegraphics[width=1\textwidth]{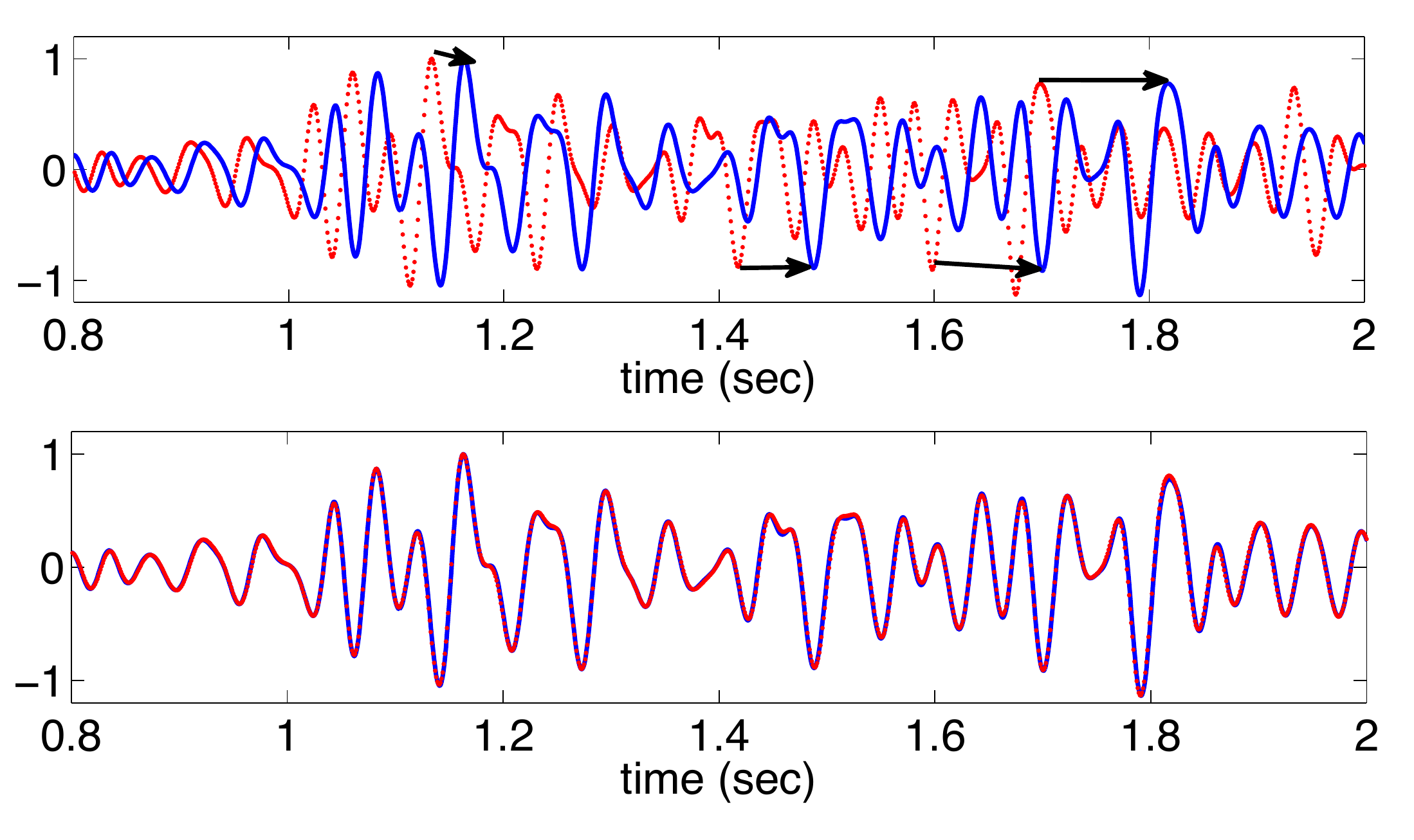} %%Marmousi_matching_noiseless.pdf}
\end{center}
\caption{
Example 1: A synthetic trace is generated using the Marmousi velocity model
and the other trace is created by warping with a known mapping. }
\label{fig:marmousi_matching}
\end{figure}

Here, we demonstrate the registration capability of the non-convex formulation (\ref{eq:uh}).
Our first example shows the registration of two {\bf noiseless} synthetic traces containing many reflected waves.
One of the two traces is obtained from a numerical experiment with the Marmousi velocity model. 
The other trace is the result of applying a warping map 
$ \displaystyle p : t \mapsto t + 0.15 \exp \left( -8 \left( t / T_c - 1 \right)^2  \right)$,
where $T_c$ is half the recording time.
Unless otherwise stated, registration is performed from zero frequency with the LFA signal $U_h$ in (\ref{eq:uh}).  
Figure \ref{fig:marmousi_matching} (bottom) shows a good registration match.

\begin{figure}
\begin{center}
\includegraphics[width=1\textwidth]{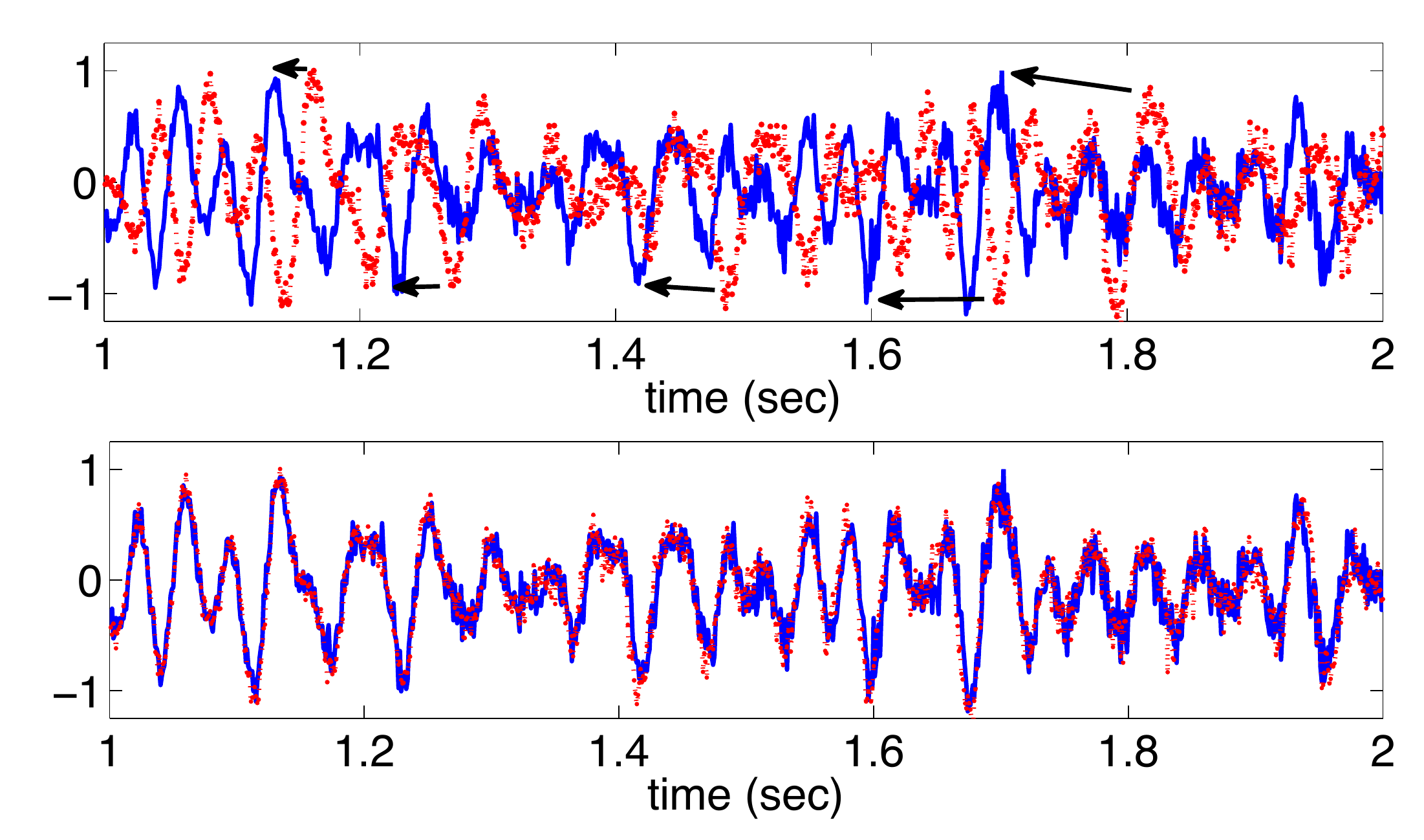}%% Marmousi_matching_noisy.pdf}
\end{center}
\caption{
Example 2: matching noisy synthetic traces. A synthetic trace is generated 
using the Marmousi velocity model
and the other trace is created by applying a known mapping.
The black arrows in the top figure connect corresponding peaks.}
\label{fig:nosiy_marmousi_matching}
\end{figure}

For the second example, we test seismogram registration with {\bf noisy} synthetic data.
A synthetic trace is obtained from  a numerical experiment with the Marmousi model and is used as reference data.
The trace is then transported using the same function as used in example 1 in order to generate predicted data.
Two independent realizations of Gaussian white noise with mean 0 and standard deviation 0.075 
are respectively added to the two traces. 
Figure \ref{fig:nosiy_marmousi_matching} shows excellent registration results in the presence of strong noise.

\begin{figure}
\begin{center}
\includegraphics[width=1\textwidth]{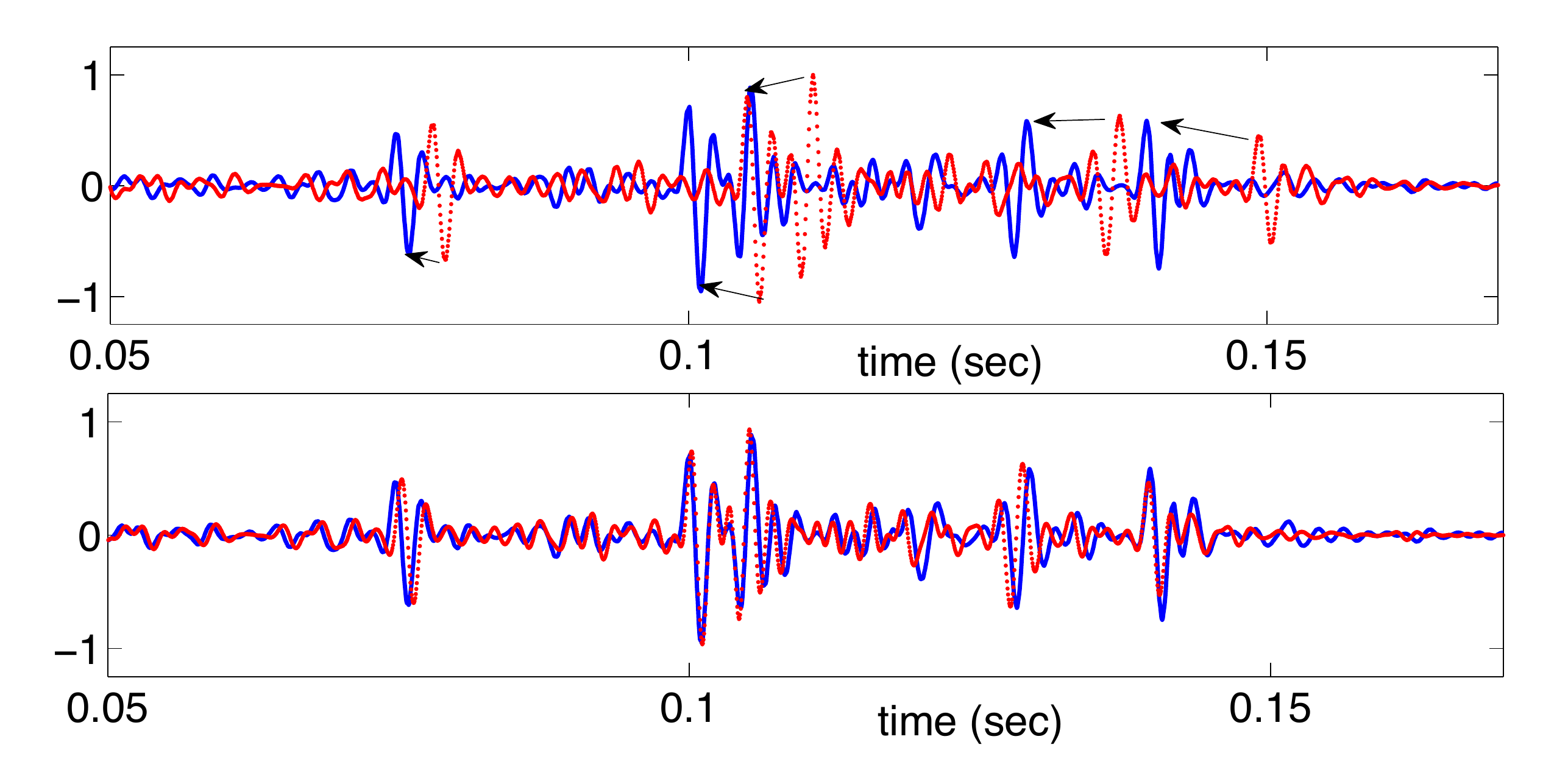} %% Marmousi_linear_transform3.pdf}
\end{center}
\caption{
Example 3: matching traces with different amplitudes and phases from the Marmousi model and a modified Marmousi model: 
(top) given observed and predicted data, 
(bottom) two traces after registration.
The black arrows indicate corresponding waves.
The predicted trace is transported towards the observed trace and its amplitude is decreased.
}
\label{fig:matching_example3}
\end{figure}

The third example uses two traces obtained from numerical experiments with {\bf two 
distinct velocity models}. An observed trace is obtained from the Marmousi velocity model $V_{Marmousi}(x,z)$ 
while we use a different velocity model $V_{pred}(x,z) = V_{Marmousi}(x,z) - 0.15 z$ for the predicted trace.
Due to the reduction in velocity, the predicted data lag behind the observed data 
up to a few wave periods in the coda.
A good, though not perfect agreement is observed between the observed trace
and the transported trace as shown in Figure \ref{fig:matching_example3}.

\begin{figure}
\begin{center}
\includegraphics[width=1.\textwidth]{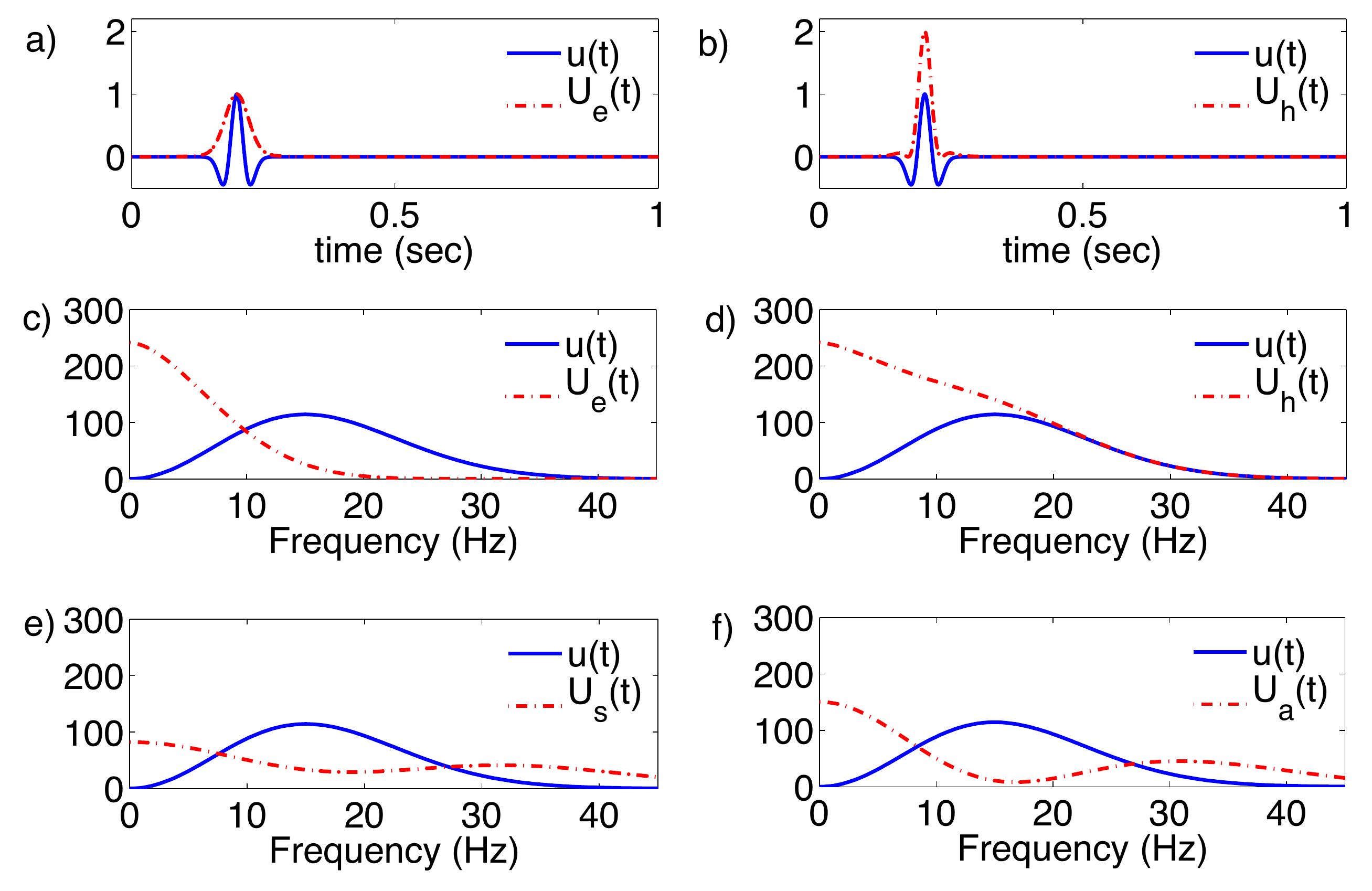} %% CreateLowFreq_wavelet_final.pdf}
\end{center}
\caption{
Spectra of signals transformed via LFA transformations (\ref{eq:uh})-(\ref{eq:ua}).
(a) $u(t)$ and its envelope $U_e(t) =| u(t) + i {\cal H} [ u(t) ]|$, 
(b) Comparison of a wavelet $u(t)$ and its LFA transformation $U_h(t)$ 
obtained using the Hilbert transform (\ref{eq:uh}),
(c) Comparison of spectra of the wavelet $u(t)$ and 
of the envelope signal $U_e(t)$, 
(d) Comparison of spectra of the wavelet $u(t)$ and of the LFA signal $U_h(t)$, 
(e) Comparison of spectra of the wavelet $u(t)$ and of the LFA signal $U_s(t) = u^2(t)$, 
(f) Comparison of spectra of the wavelet $u(t)$ and of the LFA signal $U_a(t) = |u(t)|$.
The blue solid (red dashed) lines in (c)-(f) respectively correspond to
spectra of the wavelet $u(t)$ (of the LFA signals).}
\label{fig:effect_hilbert}
\end{figure}

\subsection{Low frequency augmentation}
The examples shown in the previous section use the LFA transformation $U_h$ in (\ref{eq:uh}), 
which exploits the Hilbert transform to create low-frequency content. 
In this section, we give the rationale for our preference for $U_h$ over 
both $U_s$ and $U_a$ defined by (\ref{eq:us}) and (\ref{eq:ua}) respectively.  
In order to compare them, we look into the spectra of those LFA signals, the registration errors, and convergence of the optimization.

We examine the effect of the different LFA transformations on the frequency content of a Ricker wavelet 
with the central frequency 15 Hz.
The envelope signal $U_e(t) = | u(t) + i {\cal H} [ u( t )]| $ has strong energy in the frequency band around zero.
However, the spectrum of the envelope signal $U_e$ loses 
frequency content above the peak frequency, as shown in Figure \ref{fig:effect_hilbert} (c).
The sum of the two signals $u$ and $U_e$ should have strong energy in both frequency bands.
Our seismogram registration examples show that there is a definite improvement when using $U_h(t) = u(t) + U_e(t)$ rather than simply
$U_e(t)$ for the LFA.
Squaring an oscillatory signal of mean zero 
generates low-frequency components around zero Hz.
However, weaker frequency components 
around the center frequency 15 Hz are observed; 
this is the consequence of convolution in the frequency domain. 
Similar patterns of low-frequency augmentation and weakened signals near the central frequency 
are found in the spectra of LFA signals obtained from LFA transformation (\ref{eq:ua}).
Hence, the LFA transformation (\ref{eq:uh}) is preferred over the other two LFA transformations.

\begin{figure}
\begin{center}
\includegraphics[width=1.\textwidth]{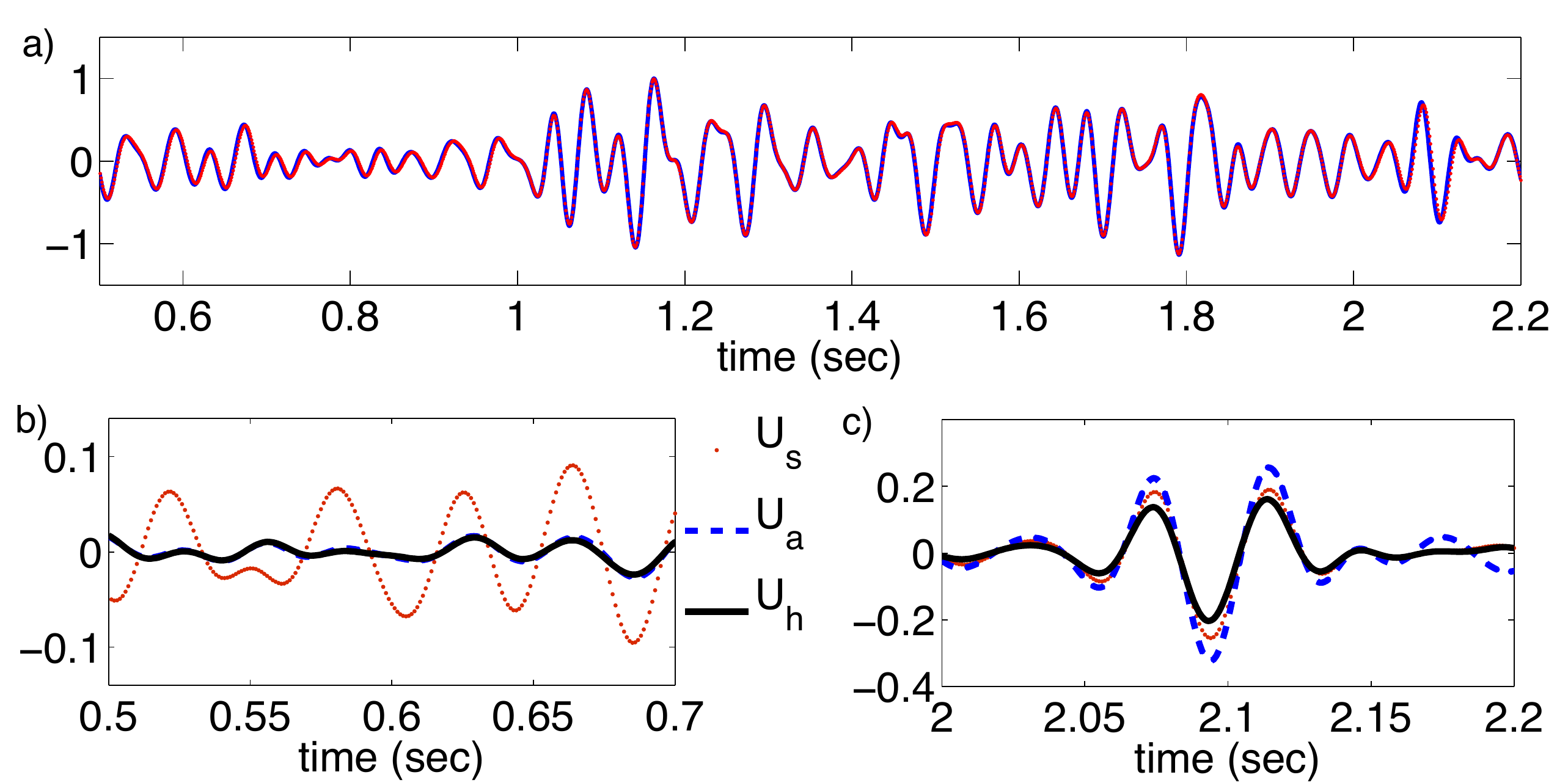} %% Error_comparison_square_hilbert_abs2.pdf}
\end{center}
\caption{
Matching synthetic data from  the Marmousi model.
The top figure shows the registration results using the LFA transform (\ref{eq:us}).
The bottom figure compares the registration errors among the three LFA signals $U_h$, $U_s$, and $U_a$.
The black solid (red marked, blue dashed) line in the bottom figures respectively corresponds to registration errors using the Hilbert
transform (squaring, absolute value).
}
\label{fig:marousi_matching_comarison}
\end{figure}

We compare the registration results obtained from the three LFA transformations (\ref{eq:uh}), (\ref{eq:us}), and (\ref{eq:ua}).
As a reference, 
the original and the transported traces using (\ref{eq:us}) are shown 
in Figure \ref{fig:marousi_matching_comarison} (top).
Registration using $U_a = |u|$ and $U_s = u^2$ yield similar results, qualitatively and performance-wise. 
The registration errors for all three LFA are plotted 
in Figure \ref{fig:marousi_matching_comarison} (bottom).
In our tests, the Hilbert-based LFA transformation enjoys the smallest registration errors among the three methods.

\begin{figure}
\begin{center}
\includegraphics[width=0.9\textwidth]{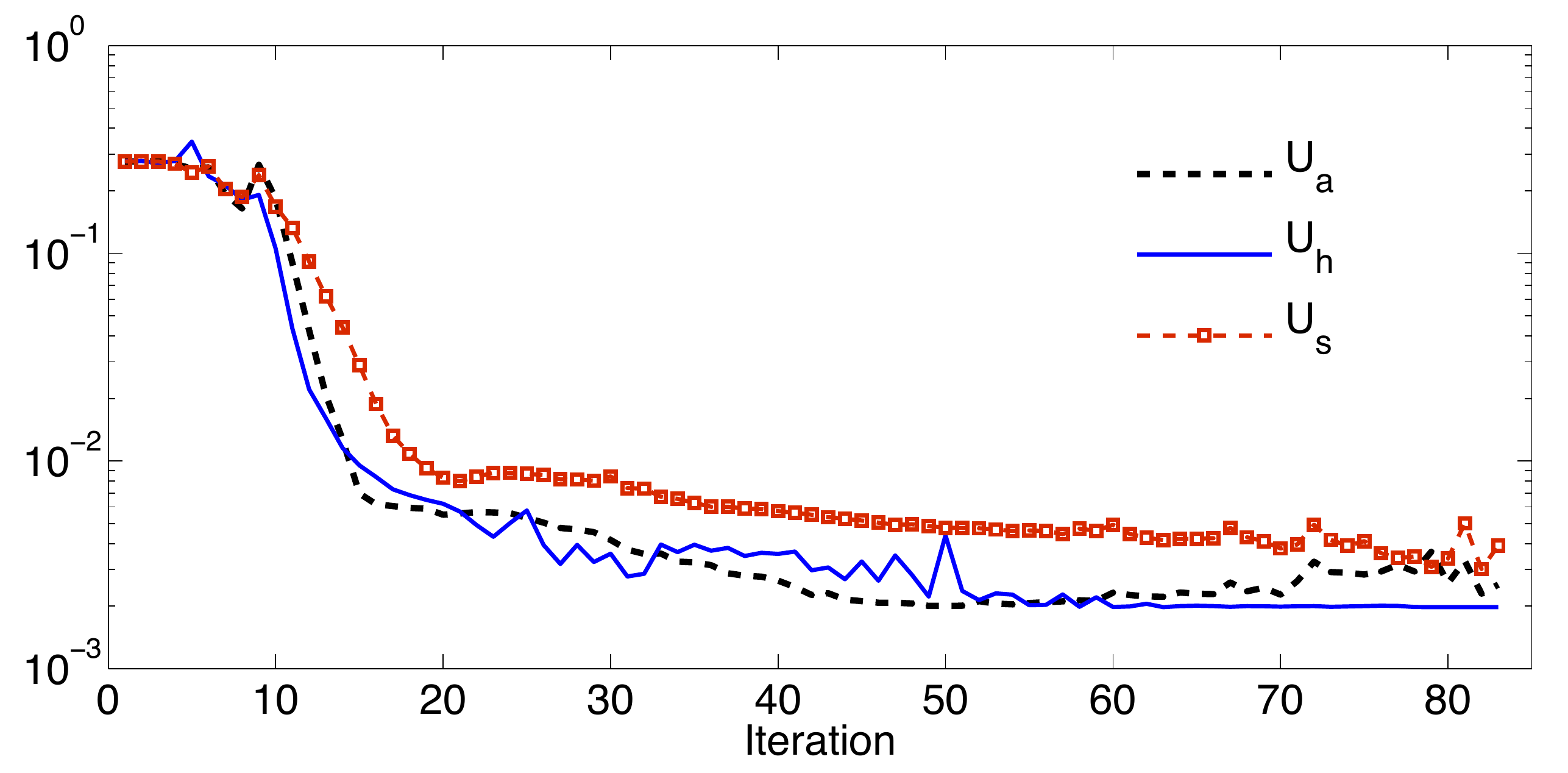} %%Cost_Comparison_3transformations3.pdf}
\end{center}
\caption{
Comparison of convergence among seismogram registration using three LFA signals: 
$U_h$, $U_s$, and $U_a$.
The black solid (red marked, blue dashed) line corresponds to 
$L^2$ data misfit using $U_h$ ($U_s$, $U_a$), respectively.
}
\label{fig:marmousi_matching_cost_decay}
\end{figure}

We also compare the decay of the registration errors among the three 
LFA transformations as registration progresses from low to high frequencies.
We observe numerically that all three LFA transformations enable seismogram registration
to converge with a reduction of the data misfit by two orders of magnitude.
The data misfits, i.e. the registration errors in the $L^2$ norm, 
for the two LFA transformations (\ref{eq:us}) and (\ref{eq:ua})
are slightly larger than those for (\ref{eq:uh}), as plotted in Figure \ref{fig:marmousi_matching_cost_decay}.

All three LFA transformations generate strong low-frequency signals, 
enabling the frequency sweep from zero frequency.
However, we claim that $U_h$ is the most appropriate among the three proposed LFA signals
in terms of frequency content, convergence of frequency sweeping, and registration errors.

\section{Numerical examples of full waveform inversion}
\label{sec:fwi_examples}

In this section, we demonstrate the potential of registration-guided least-squares (RGLS) optimization
 for waveform inversion in transmission settings with synthetic velocity models. 
Conventional LS optimization converges to a wrong velocity model in three examples with different source-receiver configurations.
On the contrary, the RGLS method can decrease the model root-mean-square (rms) error by a few orders of magnitude, in spite of temporary
increases in data misfit. We presume that the lack of monotonic decrease of the least-squares data misfit is in fact necessary for
convergence to the true model.

\subsection{Numerical methods and experimental setup}

 The acoustic scalar wave equation is used to model wave propagation in unit density media 
 with heterogeneous velocity distributions.
The equation is discretized with a
$4^{th}$ order accuracy finite difference scheme in space.
The grid size is 501 by 501 with a distance of 5 m between grid points along both directions.
Perfectly matched layers (PML) surround the computational domain. 
For the time discretization, the explicit $2^{nd}$ order leap-frog scheme is used.
A Ricker wavelet with the center frequency 50 Hz is used as an acoustic source.
This source signature is assumed to be known during inversion 
although it should be estimated in practice.
The steepest descent method is used in order to update the velocity model; 
the gradient is computed from the usual adjoint-state imaging condition involving the forward incident wavefield and the backward adjoint
field.
For RGLS optimization, 100 to 150 iterations 
suffice to significantly reduce the model rms errors.
Conventional LS optimization in many cases fails at reducing model rms errors altogether.
(We still report the result of LS optimization after a similar number of iterations.)

\begin{table}[!ht] 
\begin{center}
\begin{tabular}{l|llll} 
%\hline 
Example & Case & Reference model  &   Initial model  & Data type     \\ 
\hline \hline
\multirow{2}{*}{Example 1}& H1 & $V_{H}$    & 5100 m/s &   crosshole \\
 & L1 &   $V_{L}$     & 6000 m/s &   crosshole       \\
\hline
\multirow{2}{*}{Example 2}& H2 &   $V_{H}$    & 5100 m/s &   nearly-complete \\ 
& L2 &   $V_{L}$     & 6000 m/s &   nearly-complete \\
\hline
Example 3 & R3 &   $V_{R}$ + Noise & 5100 m/s &   nearly-complete \\  
\hline
 \end{tabular} 
\caption{Velocity models and data types for inversion examples.
Reference velocity models $V_{H}$, $V_L$, and $V_{R}$ are  
$V_{H}(x,z) = 5200 + 900 \exp( -|(x,z) -(1250,1250)|^2/10^6 )$, 
$V_L(x,z) = 5500 - 900 \exp( -|(x,z) -(1250,1250)|^2/10^6 )$, and 
$V_{R}(x,z) = 5000 + 900 \exp( -|(x,z) -(1250,1250)|^2/10^6 )$, respectively.}
\label{table:inversion_example}
\end{center} 
 \end{table}
 
The 2D velocity models with a high or low-velocity zone at the center are used in the following inversion examples 
as listed in Table~\ref{table:inversion_example}.
The velocity model in the third example contains 
noise generated by convolution of a Gaussian kernel with 
an array of normally distributed random numbers.
The initial models are constant at $V_{init}(x,z)$ = 5100 or 6000 m/s, 
i.e., without any a priori knowledge about the true models.
Such true and initial velocity models are chosen 
so that observed data do not contain caustics.
Moreover, they are intended to result in 
large discrepancies in traveltime between predicted data and observed data.

Two test configurations are used: crosshole and nearly-complete.
The data type ``crosshole" means that 
the sources are located near the left boundary edge
while the receivers are on the right side of the numerical domain.
This data type is used in example 1.
In order to separate the ill-conditioning issue from the non-convexity issue,
the second and third examples assume the availability of ``nearly-complete" data.
In this case, each shot consists of a source placed on one of the four boundary edges,
and a full set of receivers on the other three edges.

We  compare convergence of the RGLS method quantitatively
with that of the LS method for each case from H1 to R3 listed in Table \ref{table:inversion_example}.
We plot 1) true vs. converged velocity models, 2) data misfit vs. iteration count for both LS and RGLS,
and 3) rms values of $V_T - V_k$, 
where $V_T$ is the true velocity model and $V_k$ is the $k^{th}$ step velocity model.
By data misfit, we simply mean the least-squares misfit (\ref{eq:least_square_misfit}).

\subsection{Inversion example 1: crosshole tests, cases H1 and L1}
Example 1 compares the performance of RGLS and LS optimization   
using the velocity models $V_{H}$ and $V_L$ with limited data availability.
The experiment configuration for this example is similar 
to that of a crosshole test in the field; 
each shot has a source located near the left boundary ($x=10$) and 
499 equally spaced receivers near the right boundary $(x=2490,  5 \le z  \le 2495)$.
The total number of shots is 49.
Due to the limited data and aperture, it is expected that
parts of the domain are not well resolved even with the RGLS method: the inverse problem is inherently ill-posed.

\begin{figure}
\begin{center}
\includegraphics[width=0.3\textwidth]{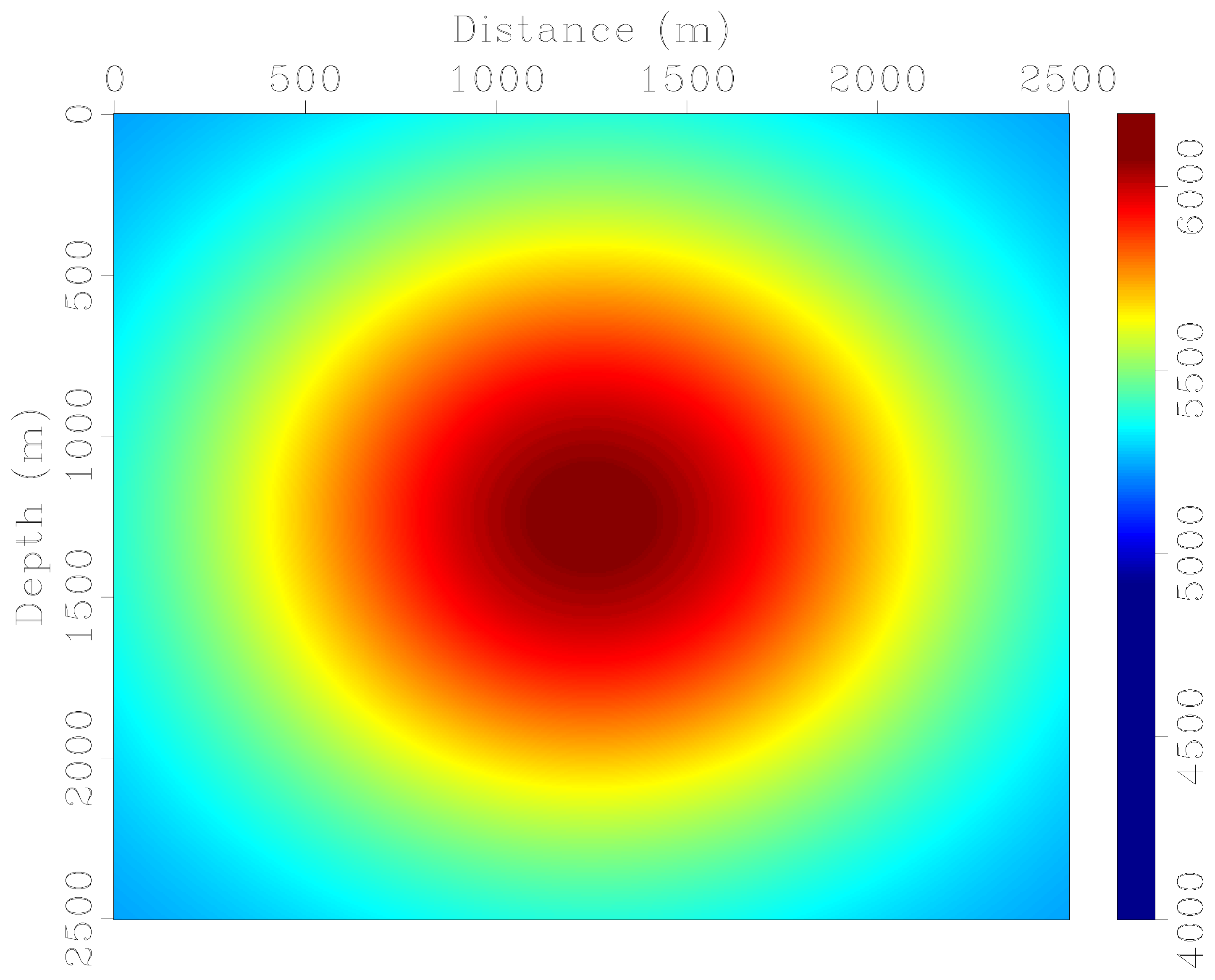} %%LPC_1side_Trans2_veltrue.pdf}
\includegraphics[width=0.3\textwidth]{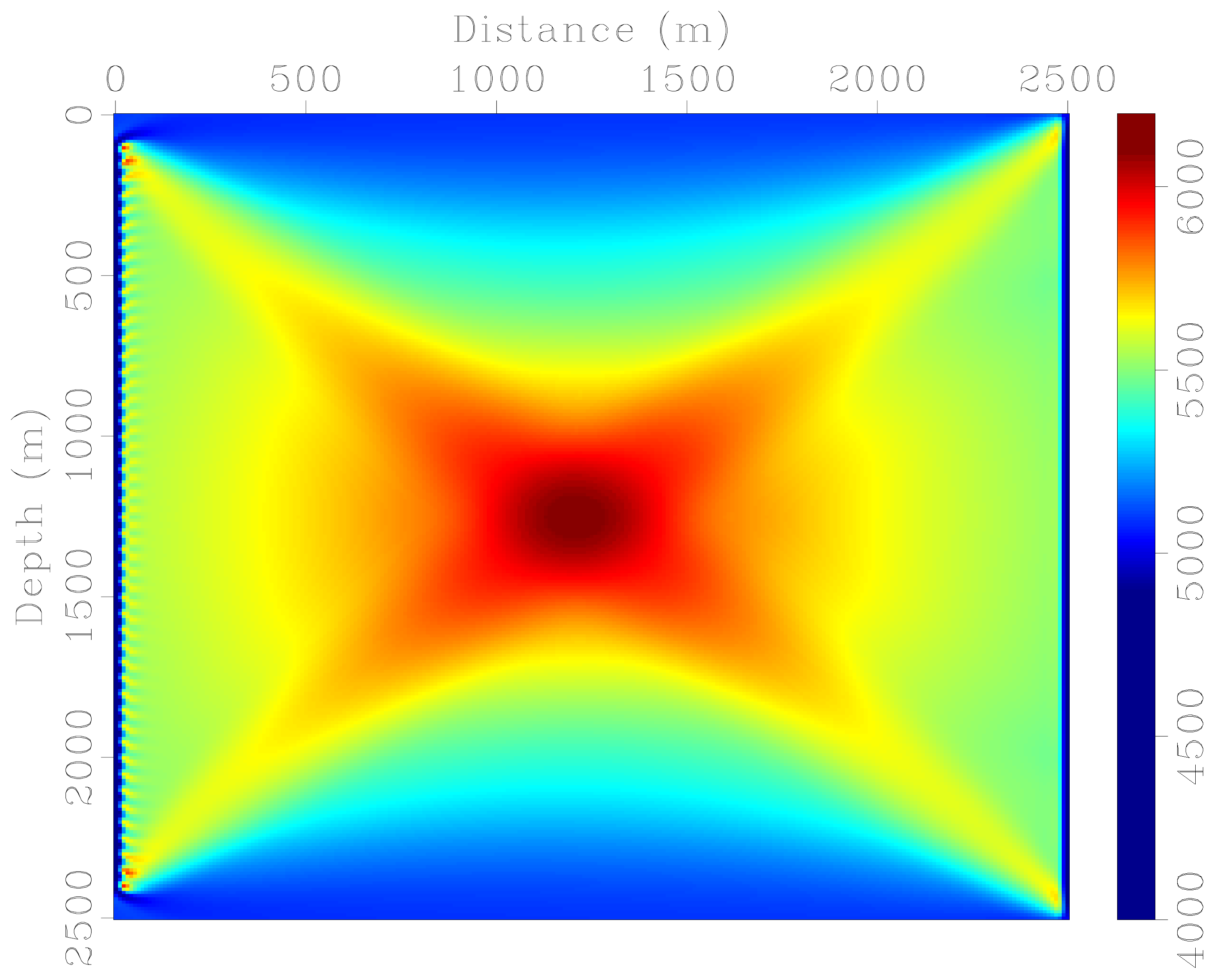} %%LPC_1side_Trans3_velsmooth130.pdf}
\includegraphics[width=0.3\textwidth]{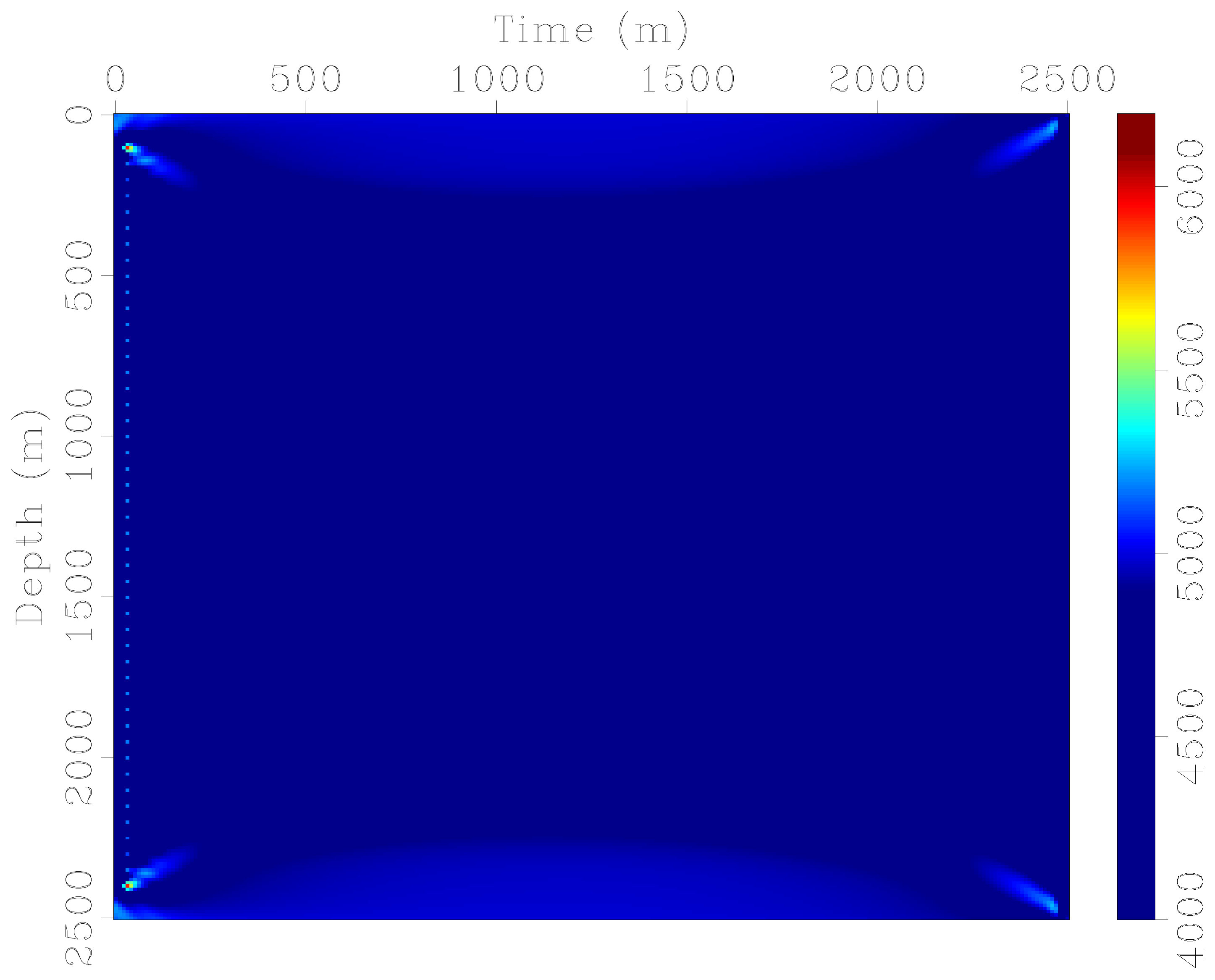} %%LPC_1side_LSQ3_velsmooth130rescale.pdf}
\includegraphics[width=0.3\textwidth]{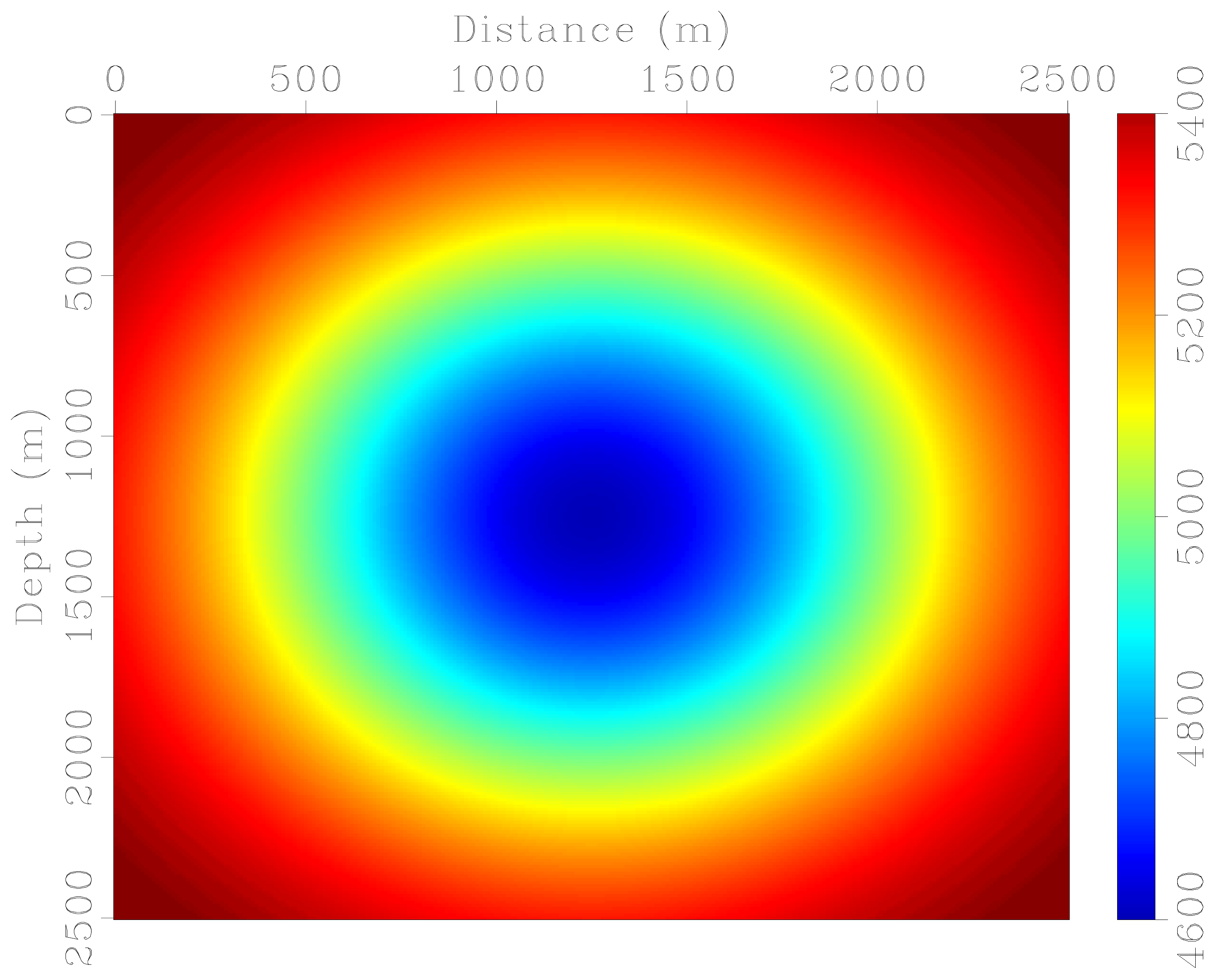} %%LMC_1side_LSQ2_veltruerescale.pdf}
\includegraphics[width=0.3\textwidth]{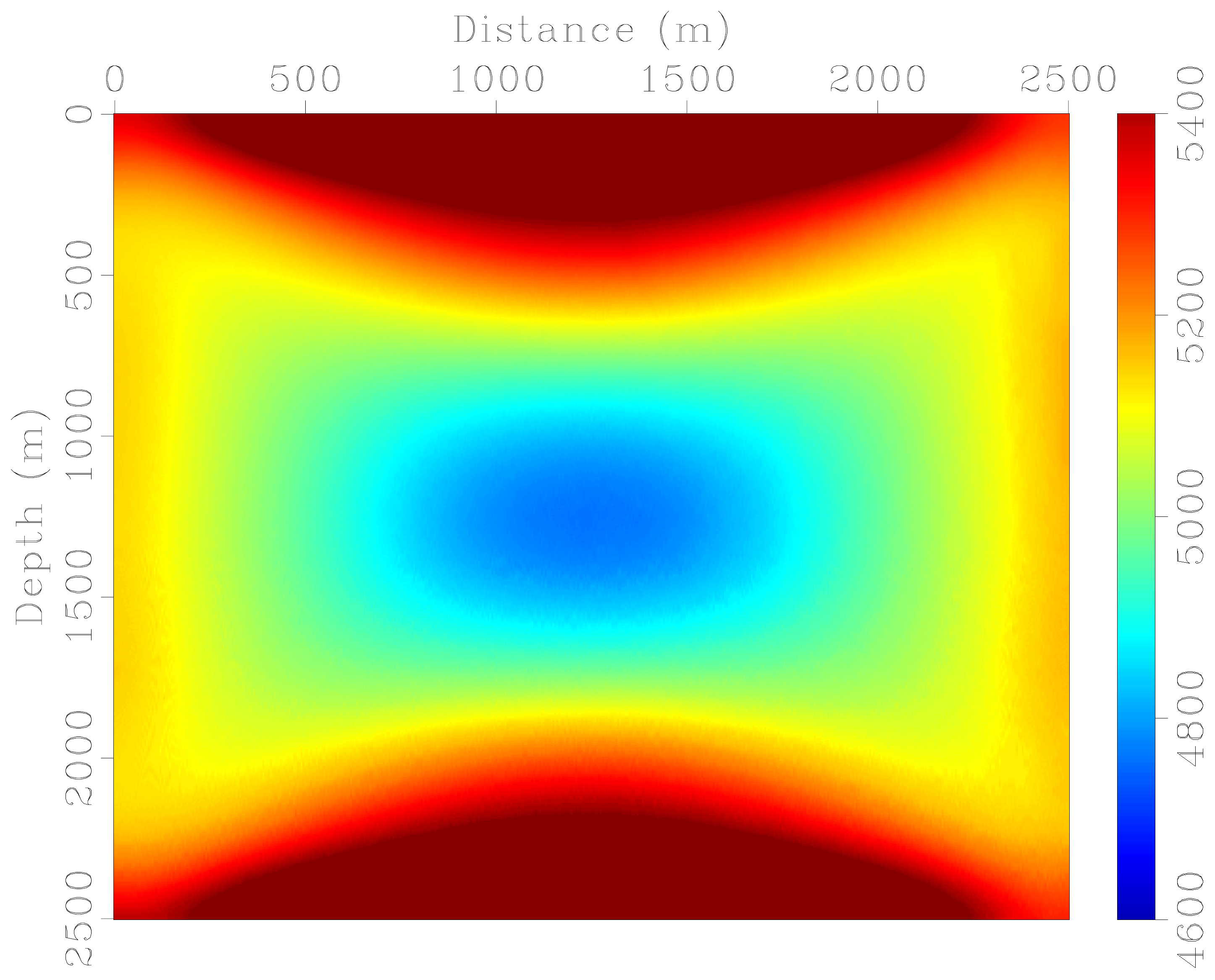} %%LMC_1side_Trans_velsmooth134rescale.pdf}
\includegraphics[width=0.3\textwidth]{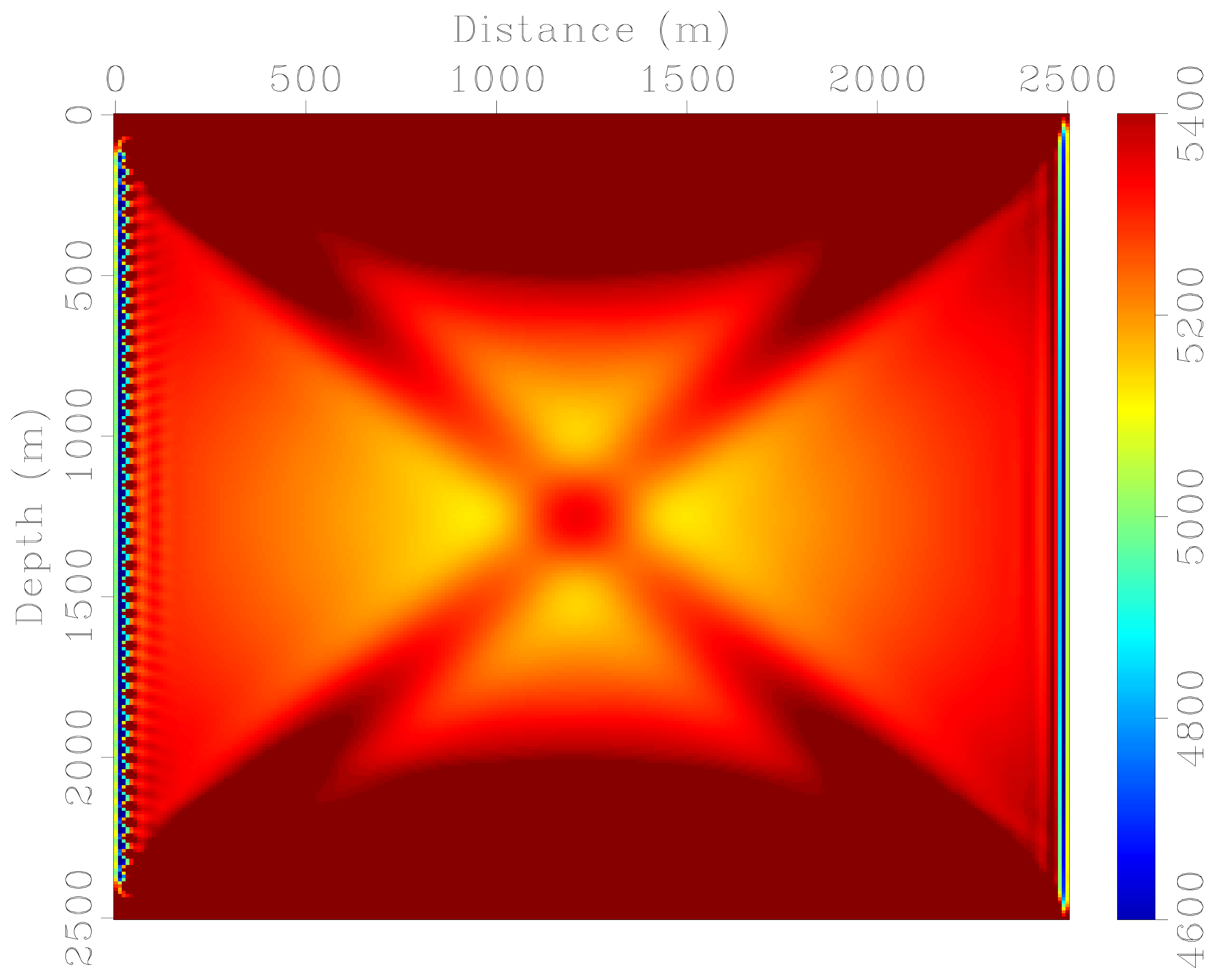} %%LMC_1side_LSQ2_velsmooth134rescale.pdf}
\end{center}
\caption{
Inversion example 1: crosshole transmission cases H1 and L1.
Plots of velocity models: (left) reference (true) models, 
(center) converged models of RGLS optimization, 
(right) converged models of LS optimization.
Top and bottom rows correspond to high and low-velocity lens models, cases H1 and L1, 
respectively.
}
\label{fig:borehole_convergence}
\end{figure}
For the H1 case (top row in Figure \ref{fig:borehole_convergence}), 
the RGLS method correctly updates the model velocity by increasing it near the center
while the LS method decreases it to below 4500 m/s.
The LS method ends up reducing the background velocity instead of increasing it. 
We understand this behavior as follows: in its drive to minimize the data misfit in the fastest possible way, 
the LS method finds it more favorable to \emph{put to zero} predicted data rather than to make it fit the observed data. 
The problem is not cycle-skipping; in this case it is simply that the time supports of the predicted and observed data are disjoint.

The two strong reflectors along the diagonals are not fully understood and deserve more investigation.
They seem to be related with the limited aperture of a crosshole test and consequent ill-conditioning.
Such reflectors do not appear in the next examples with source-receiver configurations 
involving nearly-complete data.

For the L1 case with low velocity lens reference model, 
both methods seem to work to some degree, as shown in Figure \ref{fig:borehole_convergence}, bottom row.
The converged velocity model of the RGLS method is still much closer to the reference model than
that of the LS method.
The converged model of the LS method is shown in Figure \ref{fig:borehole_convergence}, bottom right: it
is technically closer to the reference model than the initial model is, though the iterations visibly led to a spurious model (local
minimum).
The LS method suffers from the same pathological behavior of putting to zero the predicted data, as observed in the previous example.

\begin{figure}
\begin{center}
\includegraphics[width=0.47\textwidth, height=0.25\textheight]{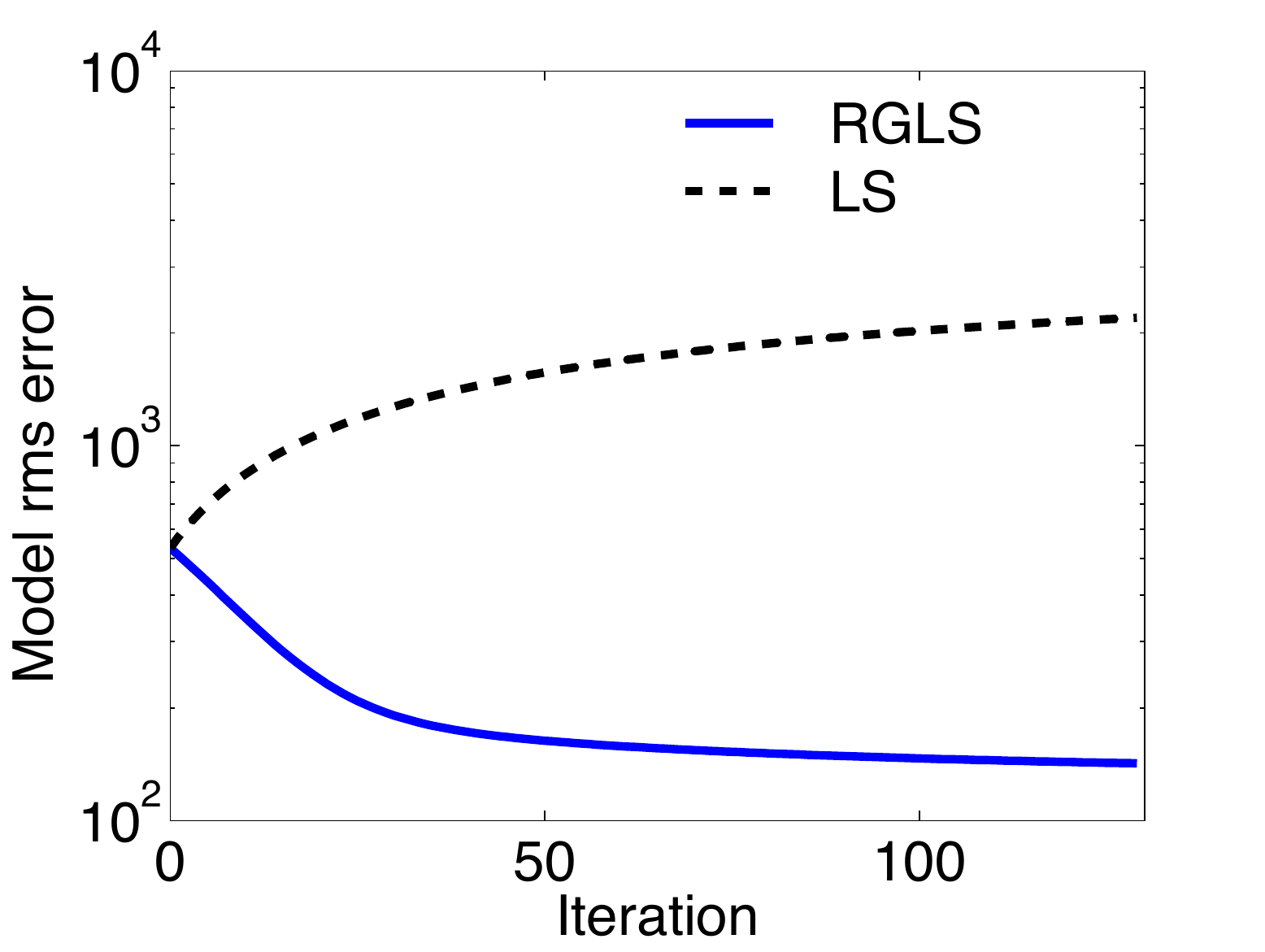} %%LPC_1side_Trans2_LSQ2_VRMSError_mod.pdf}
\includegraphics[width=0.47\textwidth, height=0.25\textheight ]{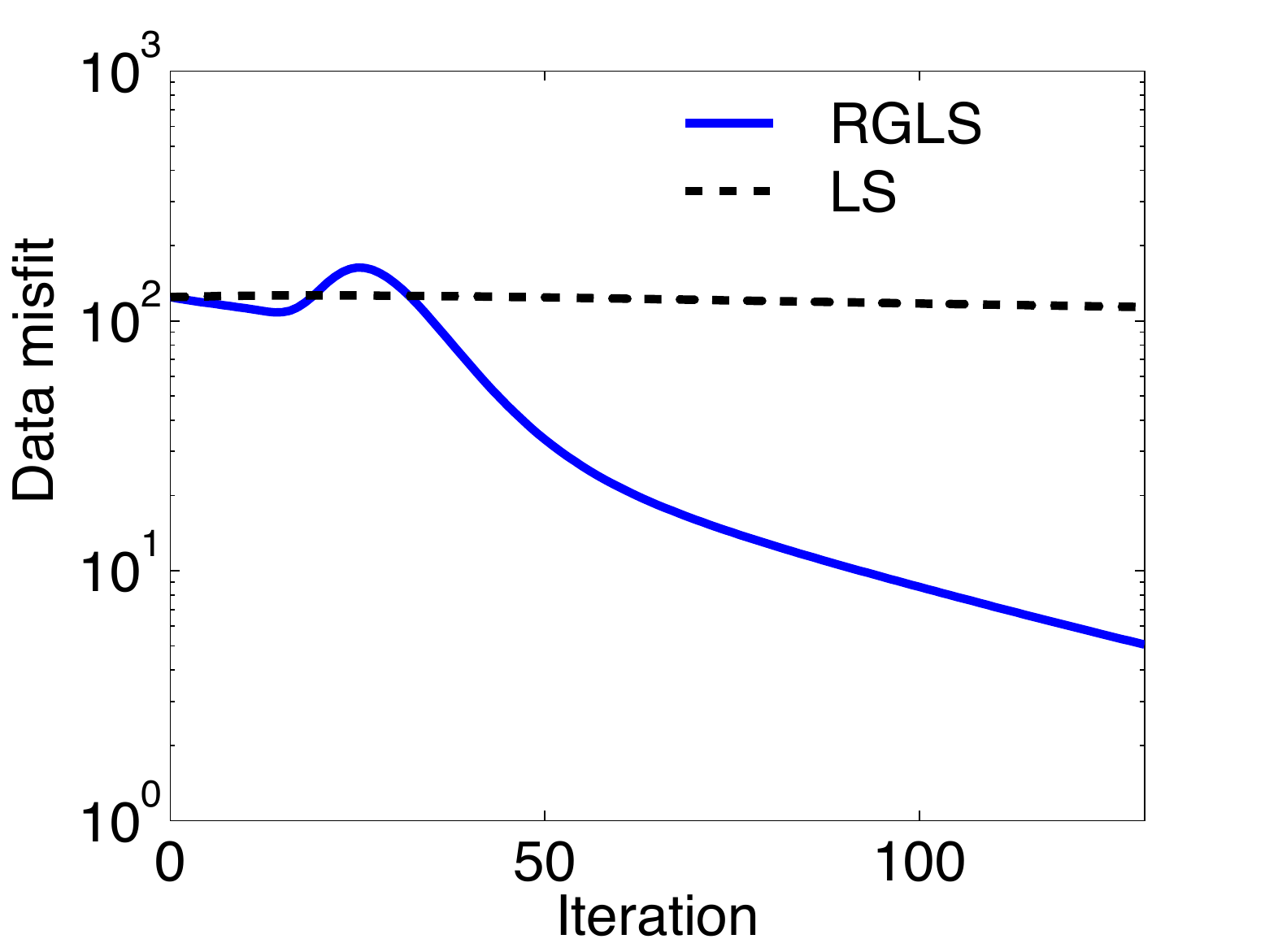} %%LPC_1side_Trans2_LSQ2_CostDecay_mod.pdf}
\includegraphics[width=0.47\textwidth, height=0.25\textheight]{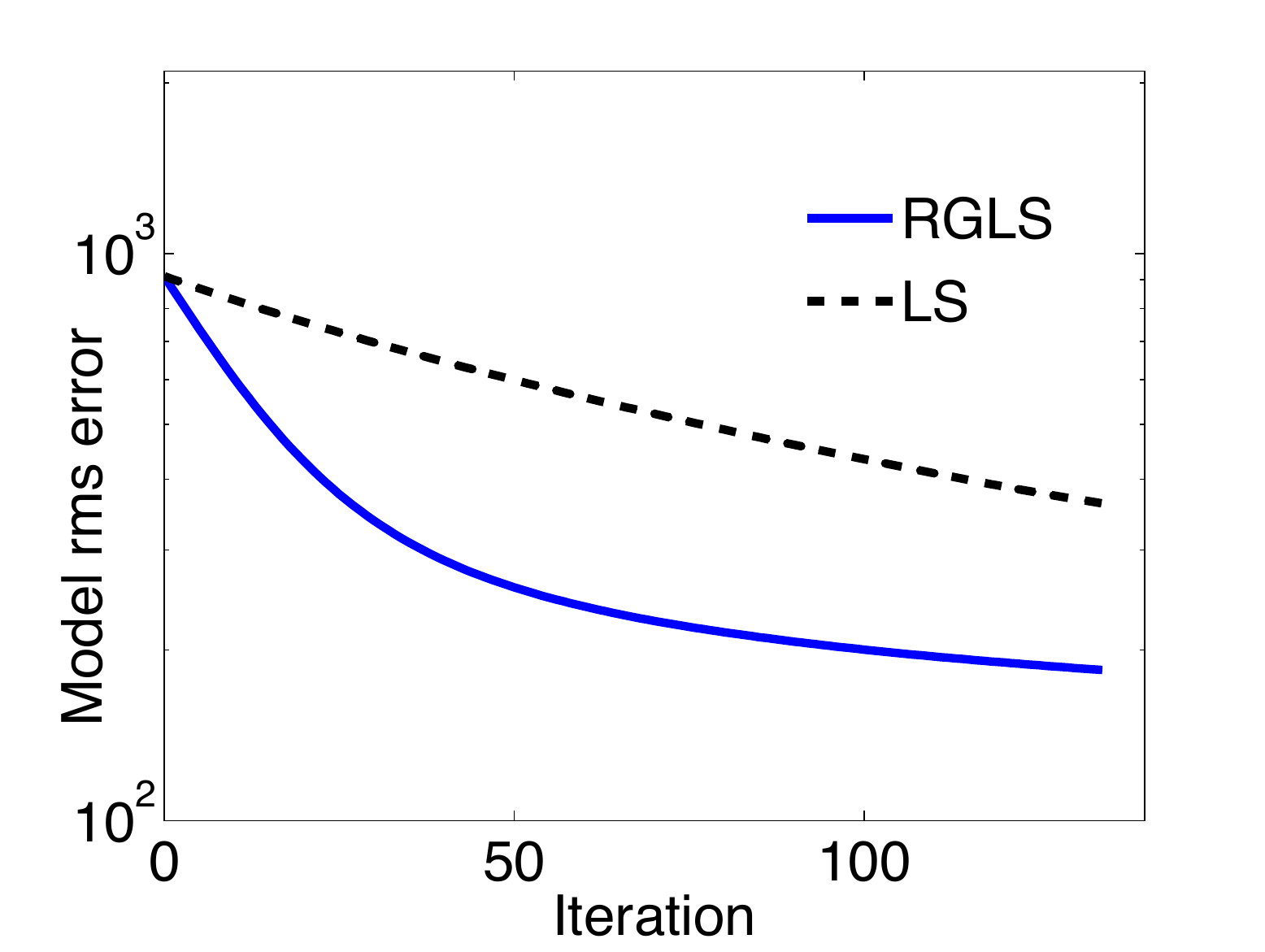} %%LMC_1side_VError.pdf}
\includegraphics[width=0.47\textwidth, height=0.25\textheight ]{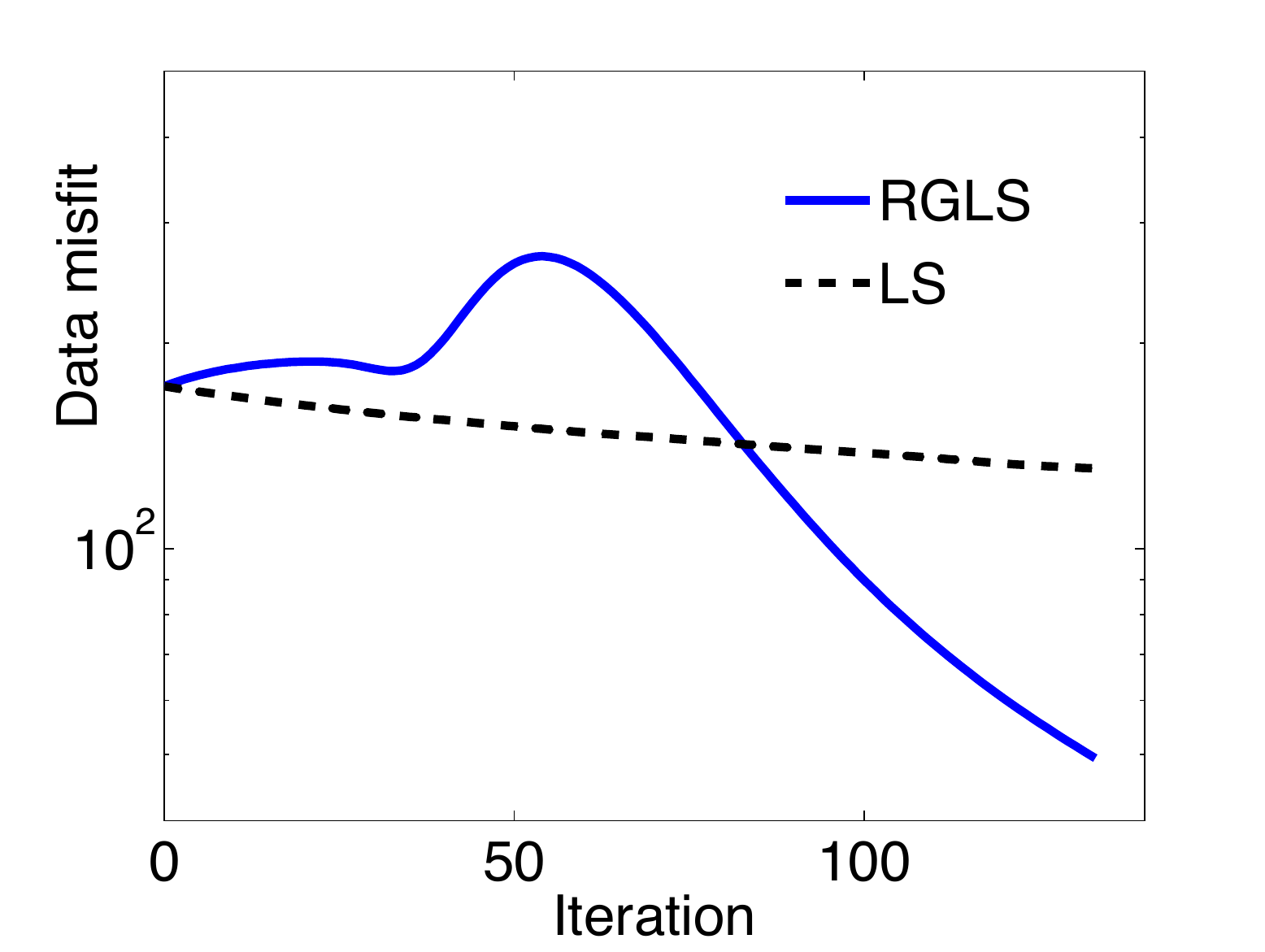} %%LMC_1side_Cost.pdf}
\end{center} 
\caption{
Inversion example 1: convergence of model rms error $V_k - V_{true} $ (left) and data misfit $J$ (right).
The top and bottom rows respectively correspond to the high-velocity H1 and low-velocity L1 lens models.}
\label{fig:borehold_failed_VelRMS_Cost}
\end{figure}

For the H1 case, LS optimization decreases the data misfit $J[m_k]$, 
but increases the model rms error gradually, 
as shown in Figure \ref{fig:borehold_failed_VelRMS_Cost} (top).
In contrast, the RGLS method temporarily increases the data misfit, but decreases the model rms error, as
the predicted data correctly move towards the observed data  during the first 20 to 30 iterations. 
For the L1 case, both RGLS and LS methods decrease model and data misfit errors, 
as shown in the bottom row of Figure \ref{fig:borehold_failed_VelRMS_Cost}. 
The RGLS optimization converges much faster than LS optimization, and converges to the correct minimum while LS gets stuck in a spurious
minimum.
  
\subsection{Inversion example 2: nearly-complete data, cases H2 and L2}

In the second example, we add more data by considering sources on all four sides of the domain. 
For any source on a given side, we consider a fine sampling of 750 receivers on the other 3 sides.
This is equivalent to a wide-aperture configuration which reduces the ill-conditioning.
The total number of shots used for computing the gradient update is 196. 
Enhancing data availability has a huge effect on the quality of converged velocity models of the RGLS method. 
The set of true and initial velocity models is unchanged from the previous inversion example.

\begin{figure}
\begin{center}
\includegraphics[width=0.3\textwidth]{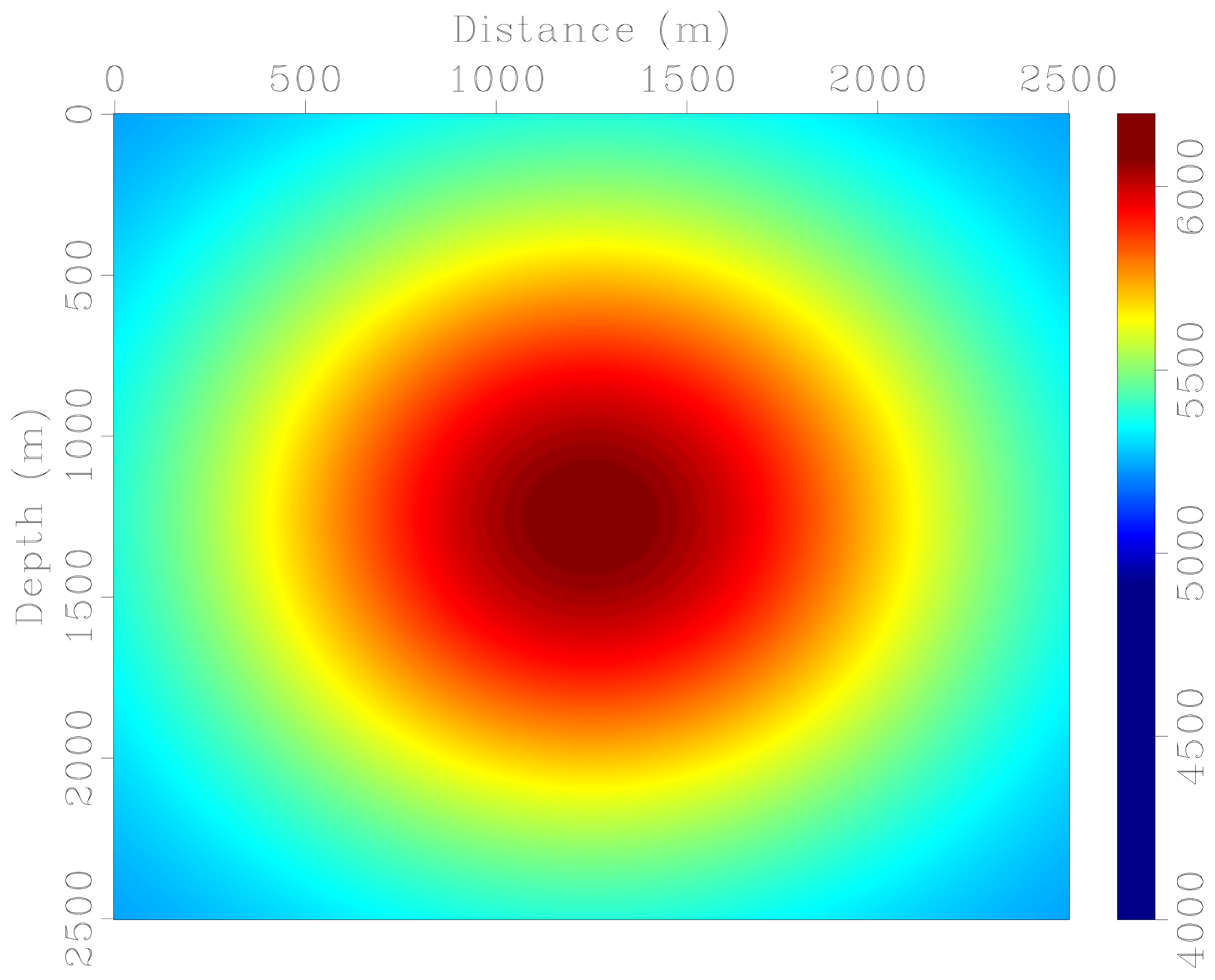} %%LPC_4sides_3sideR_Trans2_veltrue.pdf}
\includegraphics[width=0.3\textwidth]{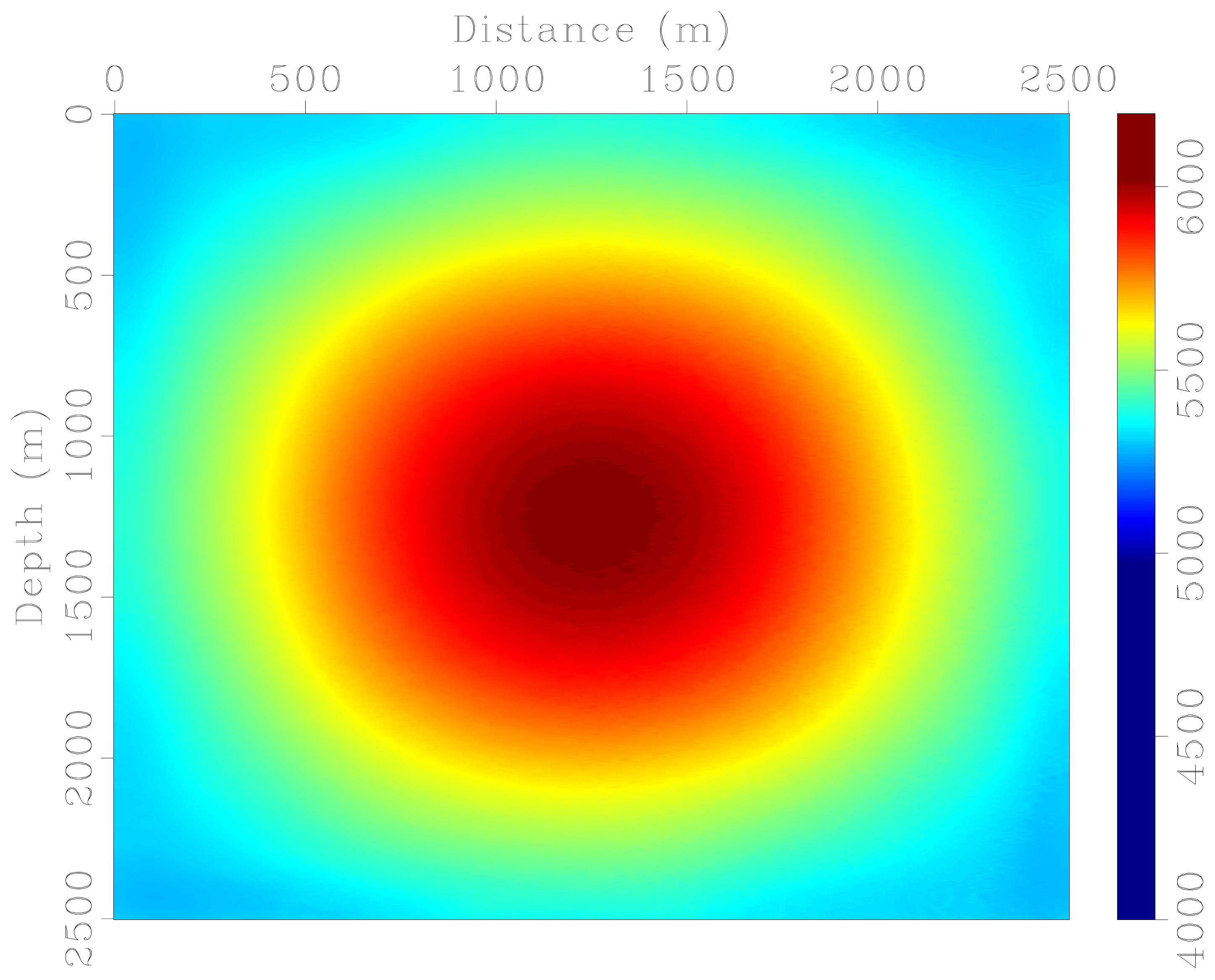} %%LPC_4sides_3sideR_Trans2_velsmooth114_redo.pdf}
\includegraphics[width=0.3\textwidth]{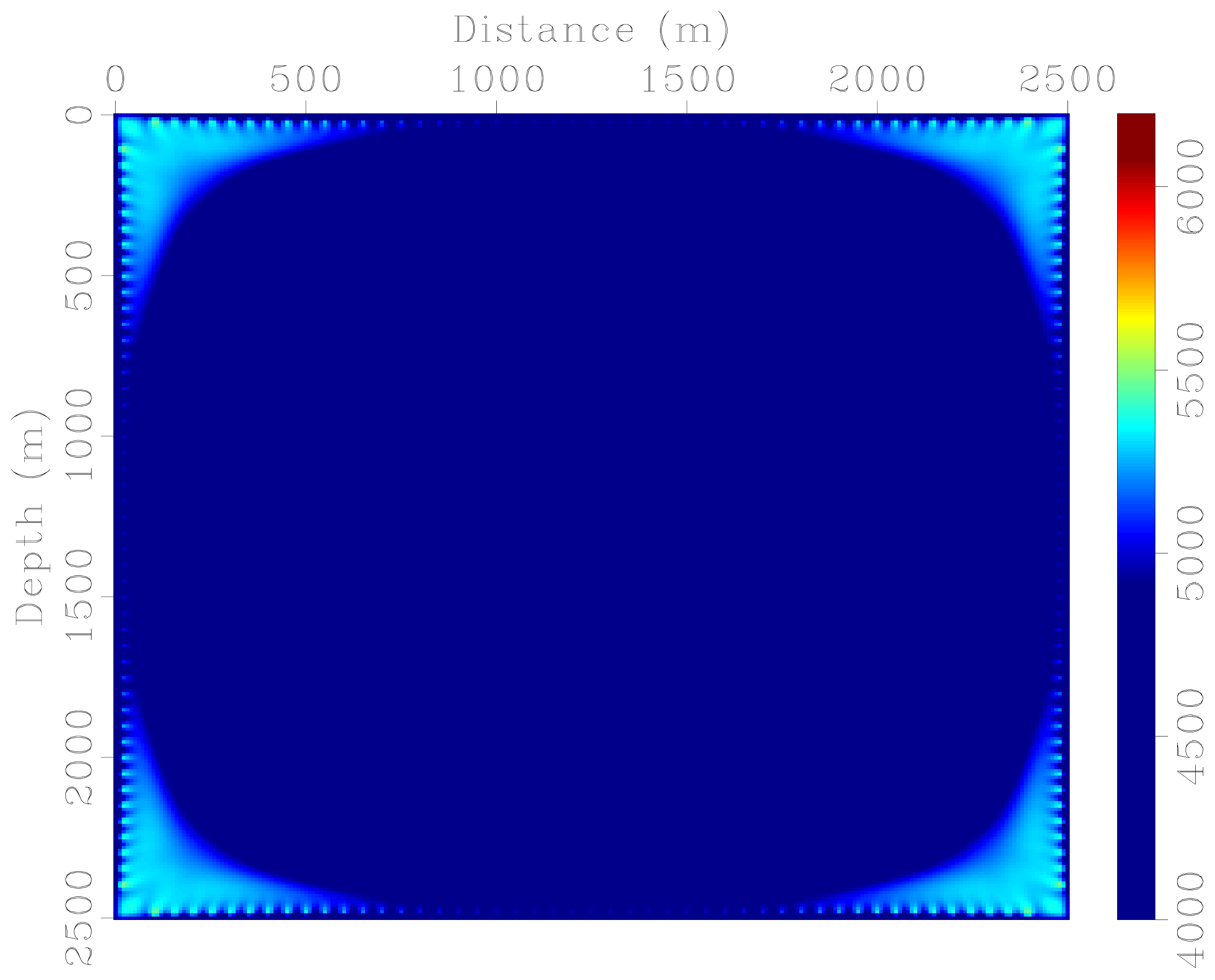} %%LPC_4sides_3sideR_LSQ3_velsmooth114rescale.pdf}
\includegraphics[width=0.3\textwidth]{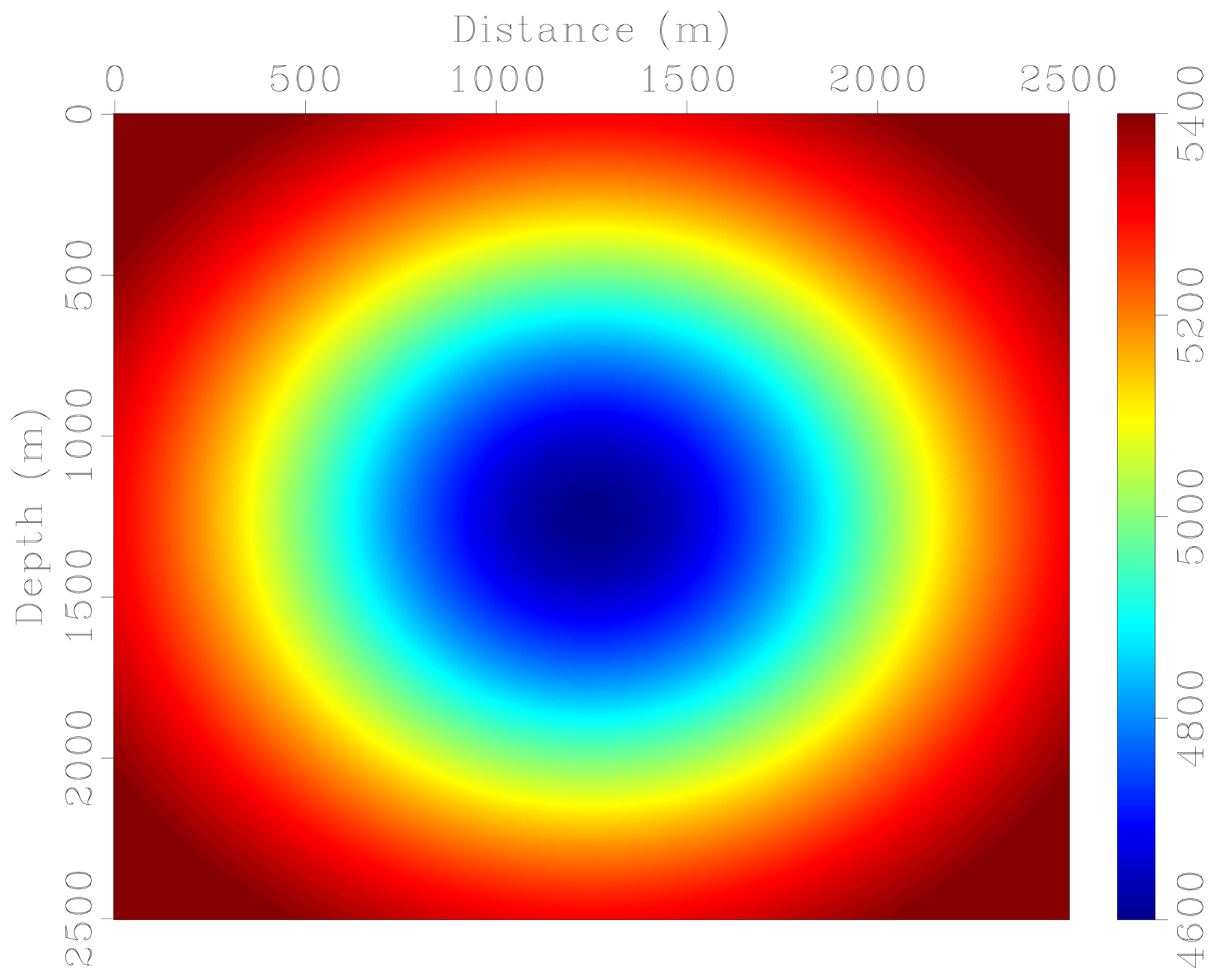} %%LMC_4sides_3sideR_Trans_veltruerescale.pdf}
\includegraphics[width=0.3\textwidth]{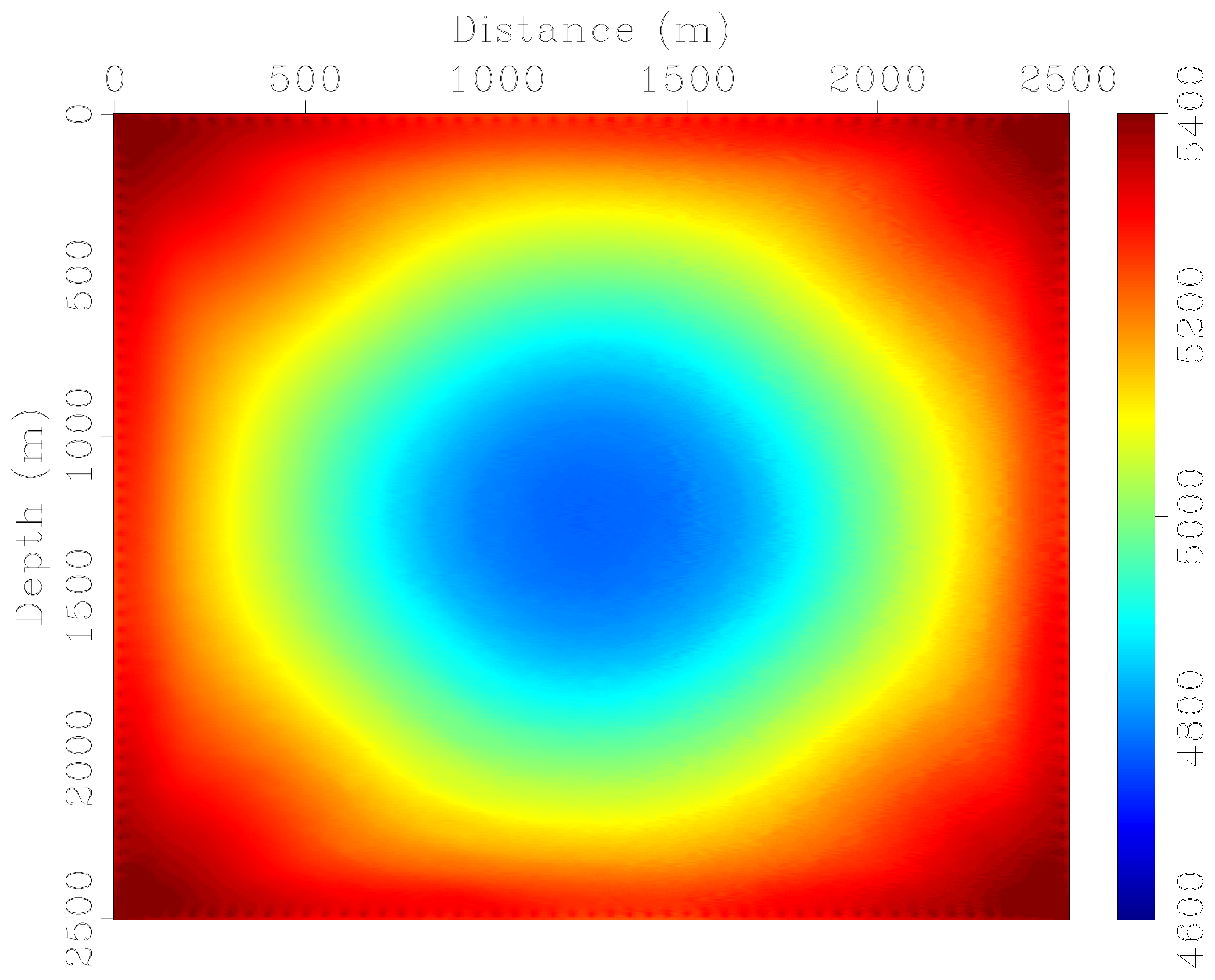} %%LMC_4sides_3sideR_Trans_velsmooth160rescale.pdf}
\includegraphics[width=0.3\textwidth]{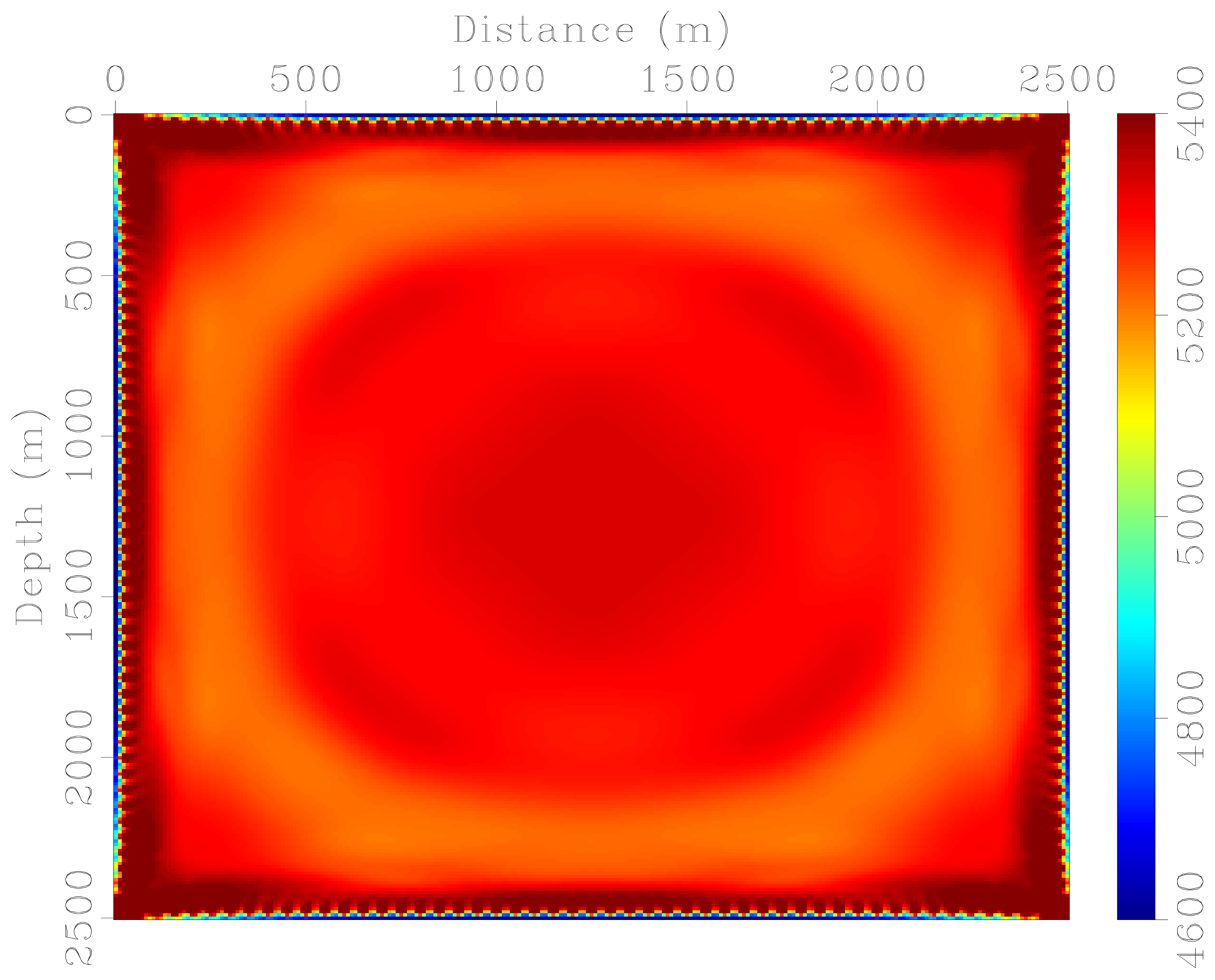} %%LMC_4sides_3sideR_LSQ2_velsmooth161rescale.pdf}
\end{center}
\caption{
Inversion example 2: plots of velocity models: (left) reference, 
(center) converged model of RGLS optimization,
(right) converged model of LS optimization. 
The top and bottom rows respectively correspond to the high-velocity and low-velocity lens models, i.e. cases H2 and L2, respectively.}
\label{fig:fully_sampled_velocity_models}
\end{figure}

For both the H2 and L2 cases, the updates of the RGLS method (100 iterations) are free of artifacts along the diagonals and the velocity
models look very close to the reference models. See Figure \ref{fig:fully_sampled_velocity_models} (center).

The LS method behaves differently in the H2 case than in the L2 case:
the method converges to a wrong model in the former case, and to the correct model in the latter case (though after a very large number of
iterations). 
For the H2 case with a high velocity zone, 
LS optimization decreases the velocity model, hence updates the velocity model in the wrong direction. 
We attribute this behavior to the fact that the LS method seeks to make the predicted waves arrive later and with smaller amplitudes, 
hence again trying to put the prediction to zero rather than matching it to observed data.
As a result, the model velocity reaches around 3500 m/s at the center, 
as shown in Figure \ref{fig:fully_sampled_velocity_models} (top right), 
which is much lower than the initial velocity 5100 m/s.  
(However, notice that the four corners seem to be updated properly. 
There is no cycle-skipping there: predicted data with short travel times are within a quarter of a wave length from corresponding observed
data around the corners. )

In the L2 case, the LS method updates correctly parts of the domain near the boundary as well as the corners,
as shown in Figure \ref{fig:fully_sampled_velocity_models} (bottom right).
Interestingly, a hundred more iterations result in a much better
velocity model with large errors only at the center. The area of the error zone at the center gets smaller 
as the iterations proceed. 
We understand this behavior as due to the fact that LS minimization whittles away at the center of the image by progressively matching data
from the corners inward.

\begin{figure}
\begin{center}
\includegraphics[width=0.45\textwidth, height=0.25\textheight]{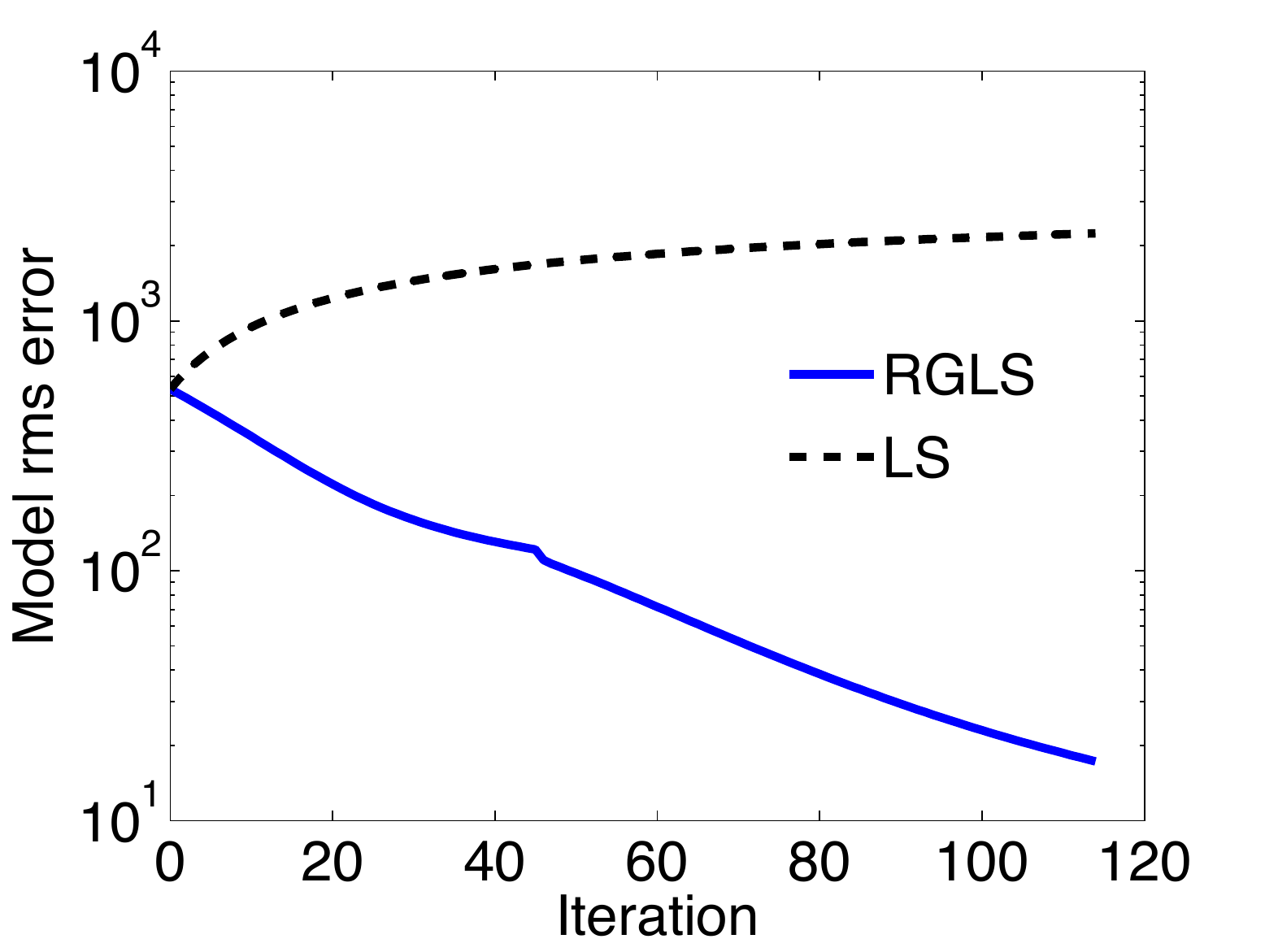} %%LPC_4sides_3sideR_VRMSError_redo.pdf}
\includegraphics[width=0.45\textwidth, height=0.25\textheight ]{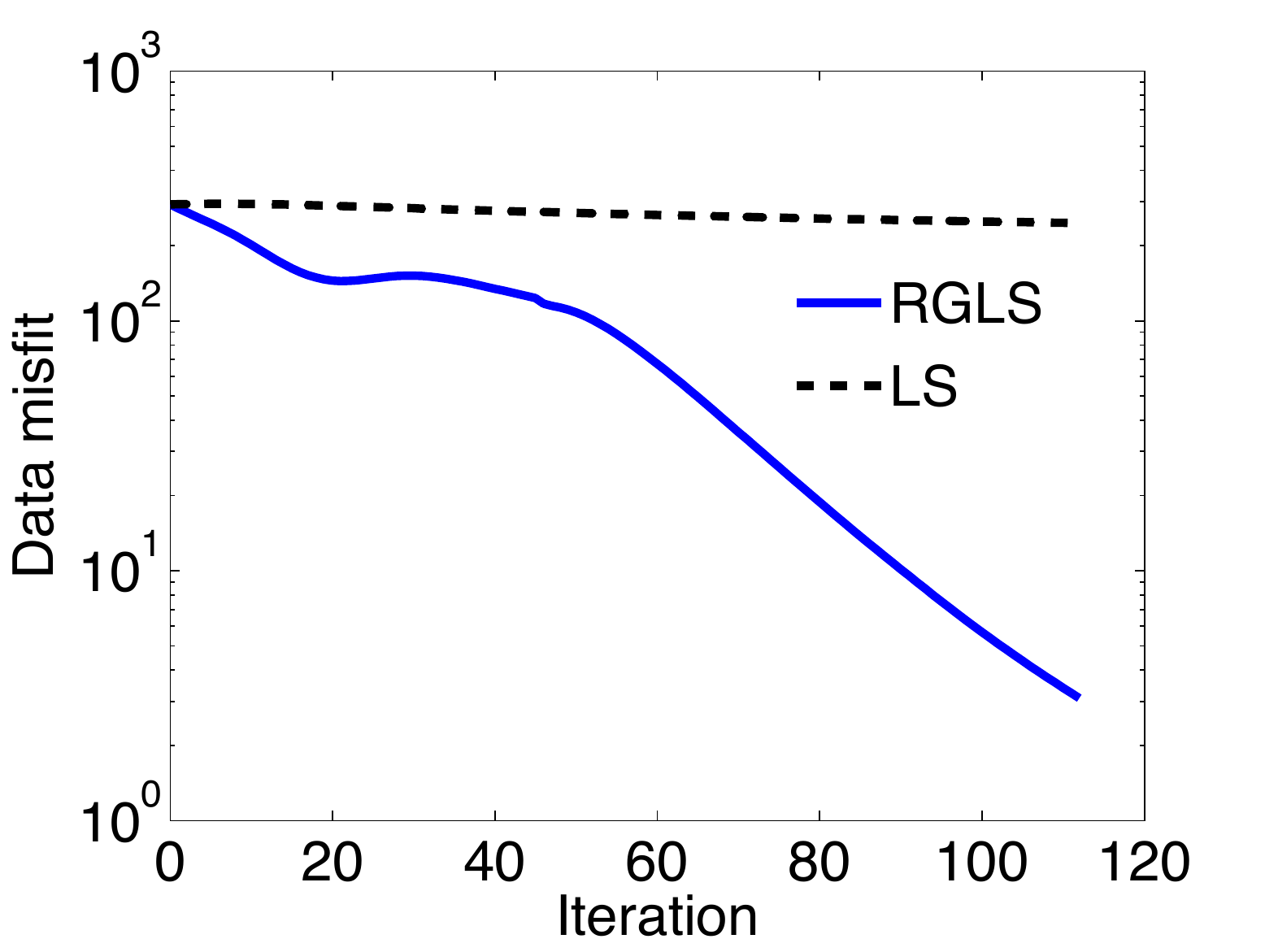} %%LPC_4sides_3sideR_CostDecay_redo.pdf}
\includegraphics[width=0.45\textwidth, height=0.25\textheight]{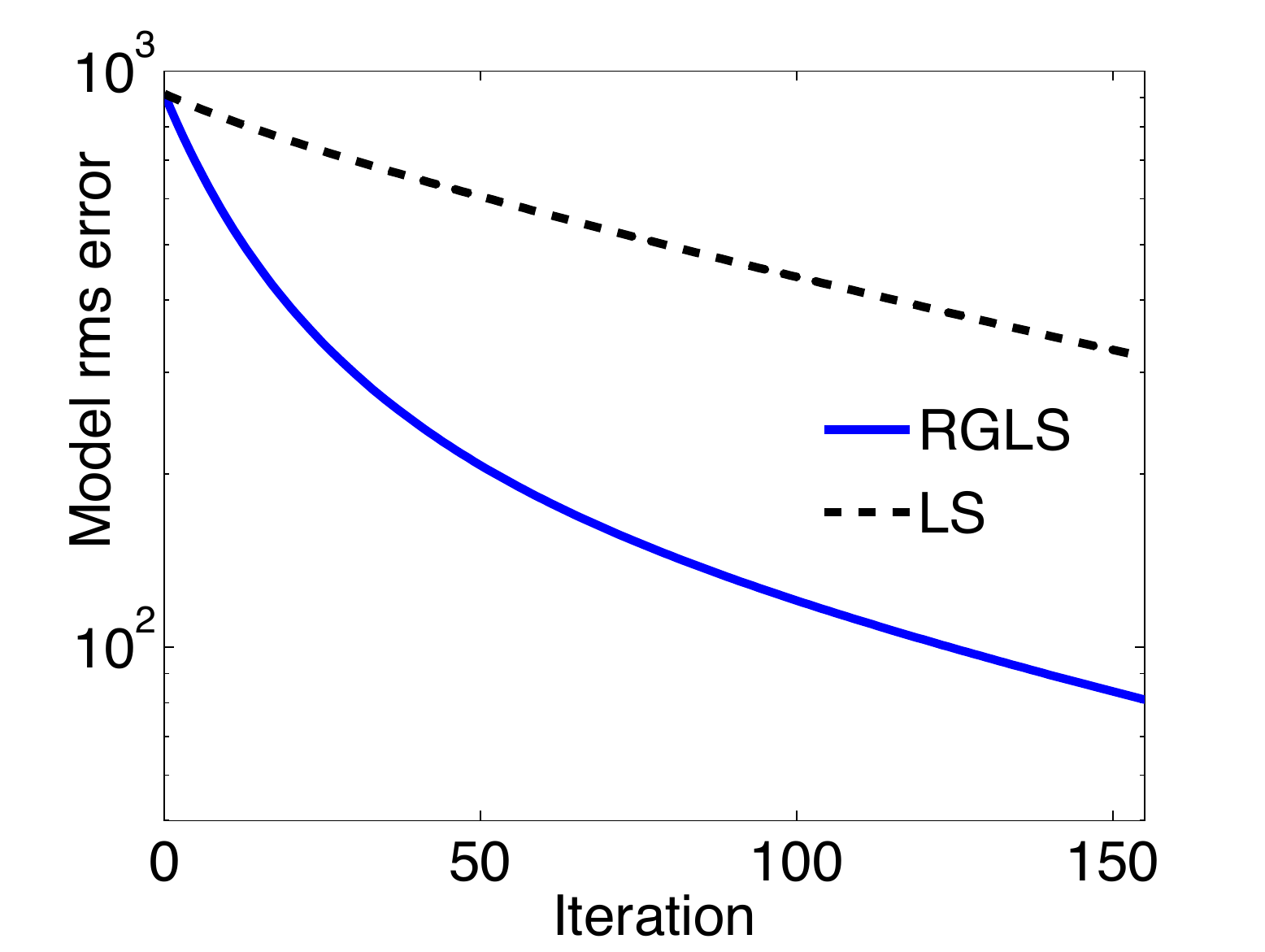} %%LMC_4sides_3sideR_VelRMS.pdf}
\includegraphics[width=0.45\textwidth, height=0.25\textheight ]{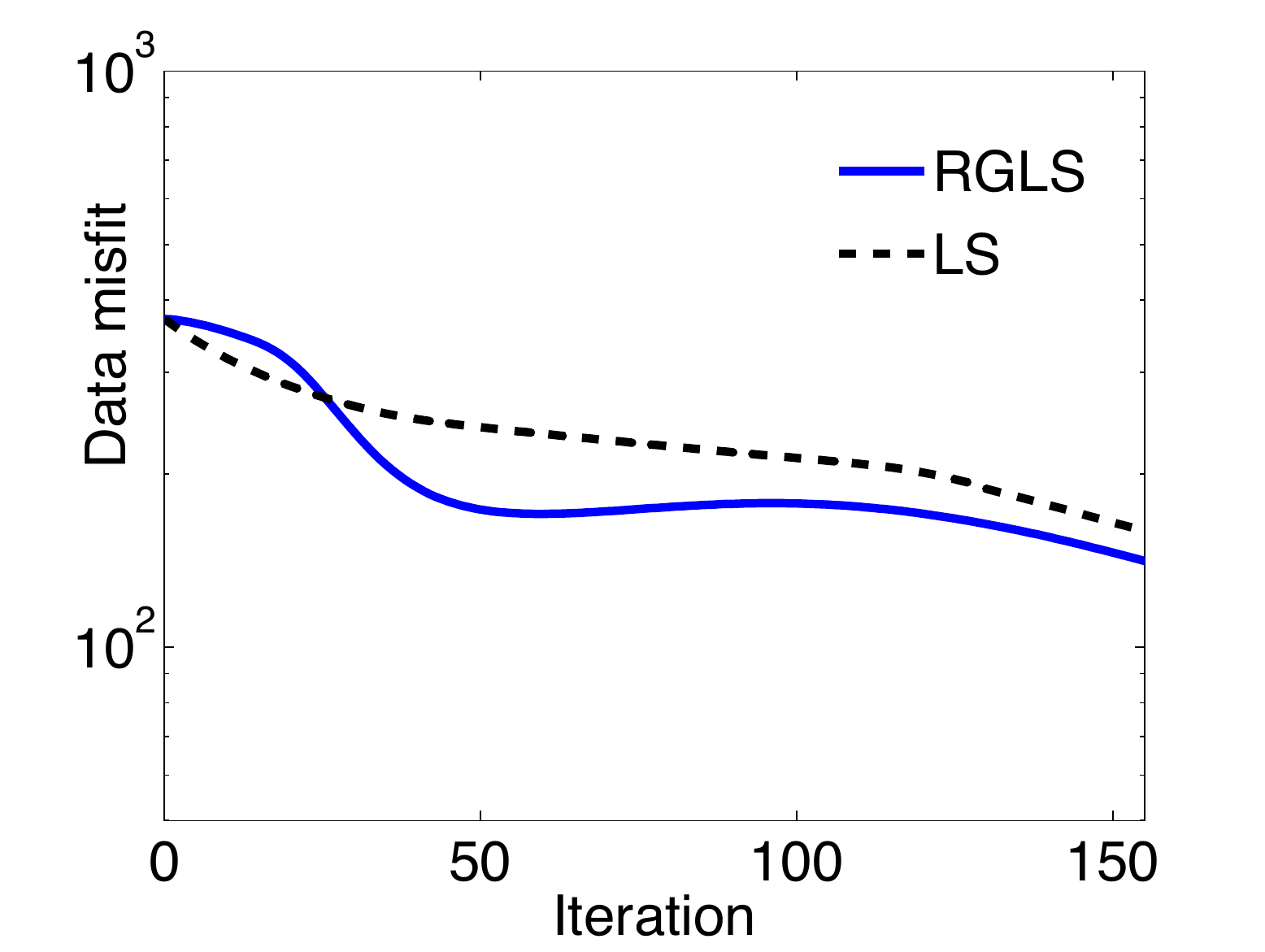} %%LMC_4sides_3sideR_Datamisfit.pdf}
\end{center} 
\caption{Inversion example 2: 
convergence of model (velocity) rms error $V_k - V_{true} $ (left) and data misfit $J$ (right)
The top and bottom rows respectively correspond to the high-velocity H2 and low-velocity L2 lens models.}
\label{fig:fully_sampled_convergence_behavior}
\end{figure}

For the LS method in the H2 case,  the data misfit decreases but the model velocity error
increases gradually, as shown in Figure \ref{fig:fully_sampled_convergence_behavior} (top).
RGLS optimization correctly updates the velocity model, decreasing the model (velocity) error, 
and ``flies above the nonconvex landscape" of the LS misfit cost function. 
As in the previous example, a hump is observed in the $L^2$ data misfit.
In the L2 case, both RGLS and LS converge, reducing both data and model errors.
However, RGLS is much faster. LS updates the model in a slow, piecemeal way from the boundary inwards.

\subsection{Inversion example 3: nearly-complete data, case R3}

\begin{figure}
\begin{center}
\includegraphics[width=0.3\textwidth]{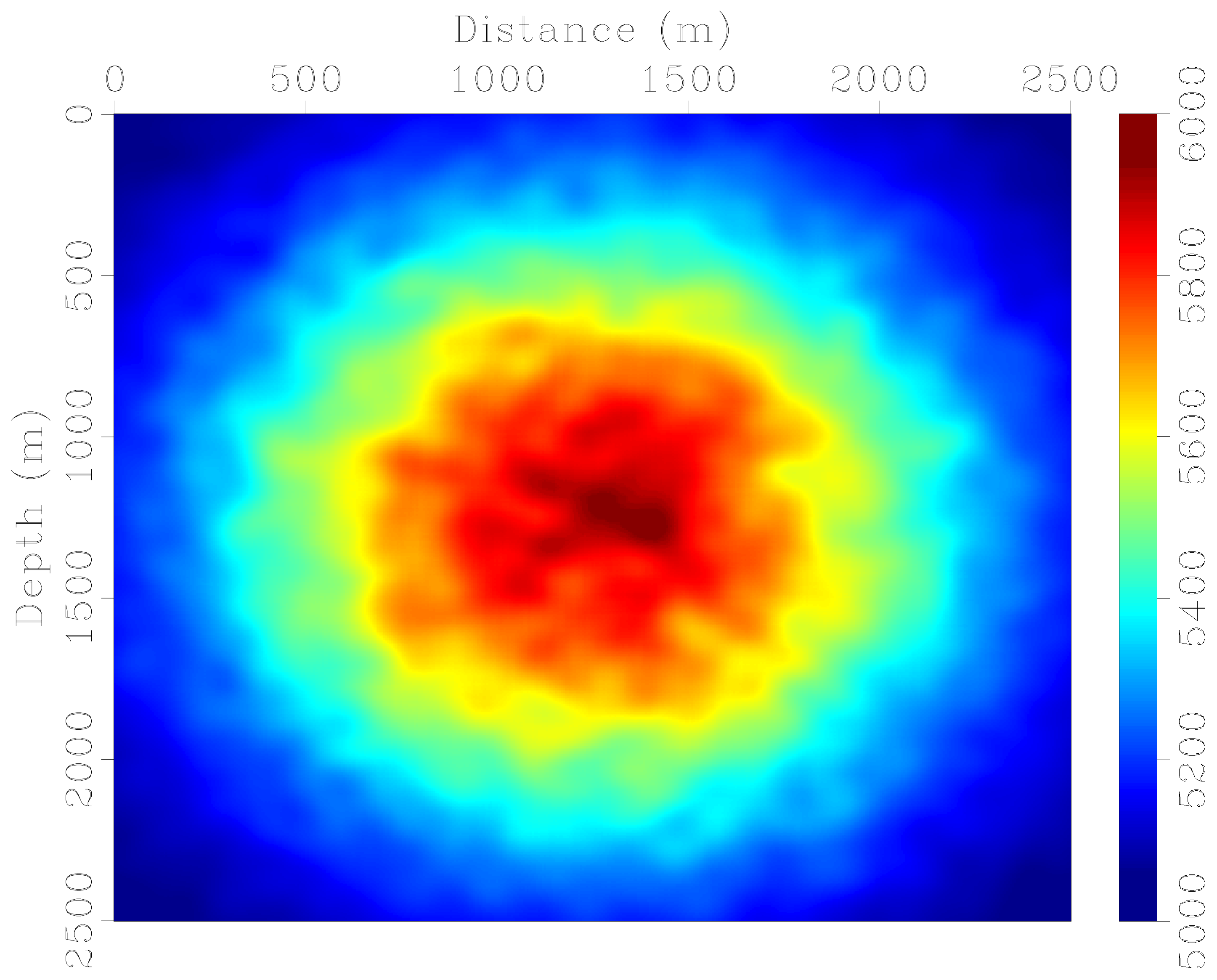} %%LPC_4s_Trans1_veltrue.pdf}
\includegraphics[width=0.3\textwidth]{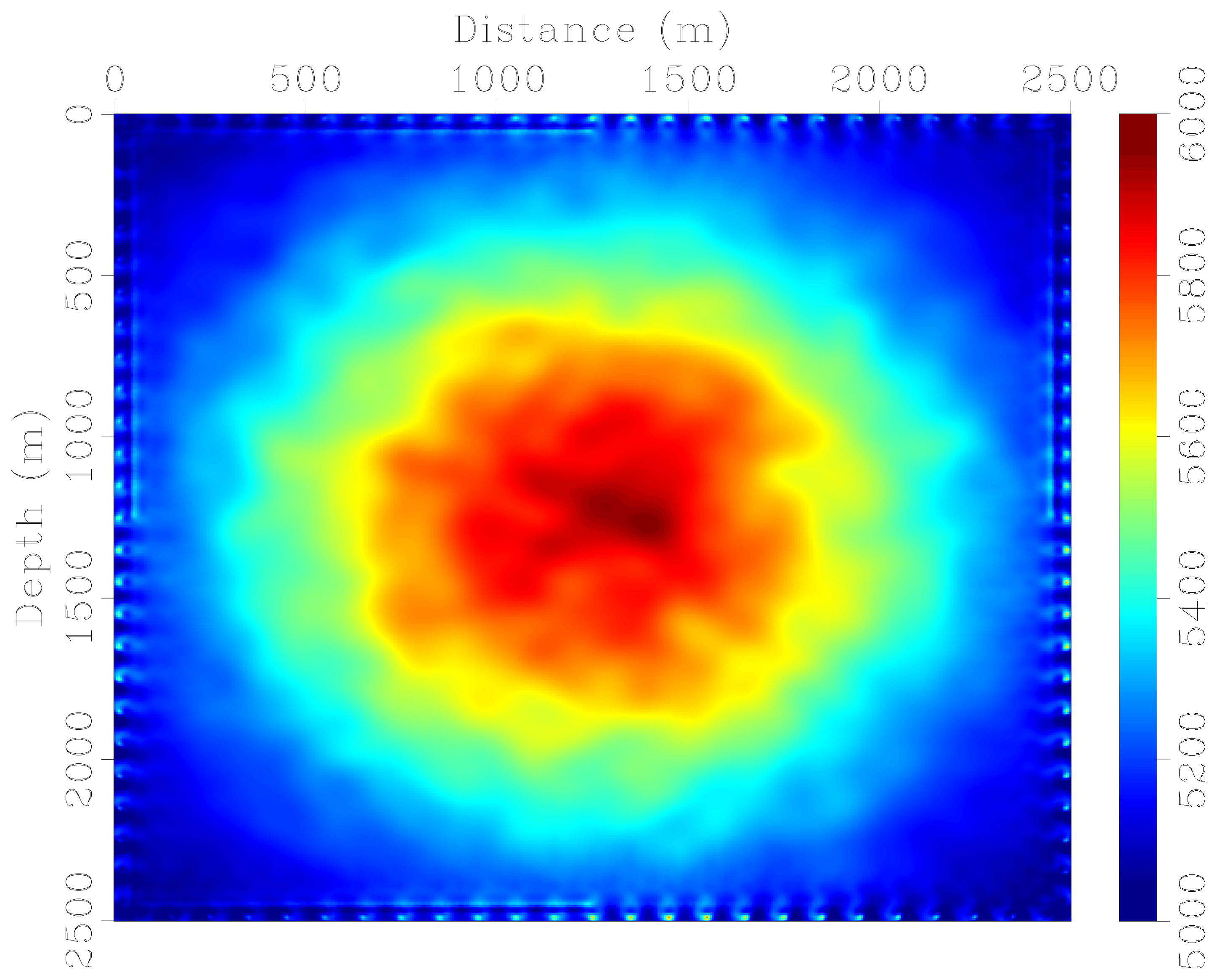} %%LPC_4s_Trans1_velsmooth100.pdf}
\includegraphics[width=0.3\textwidth]{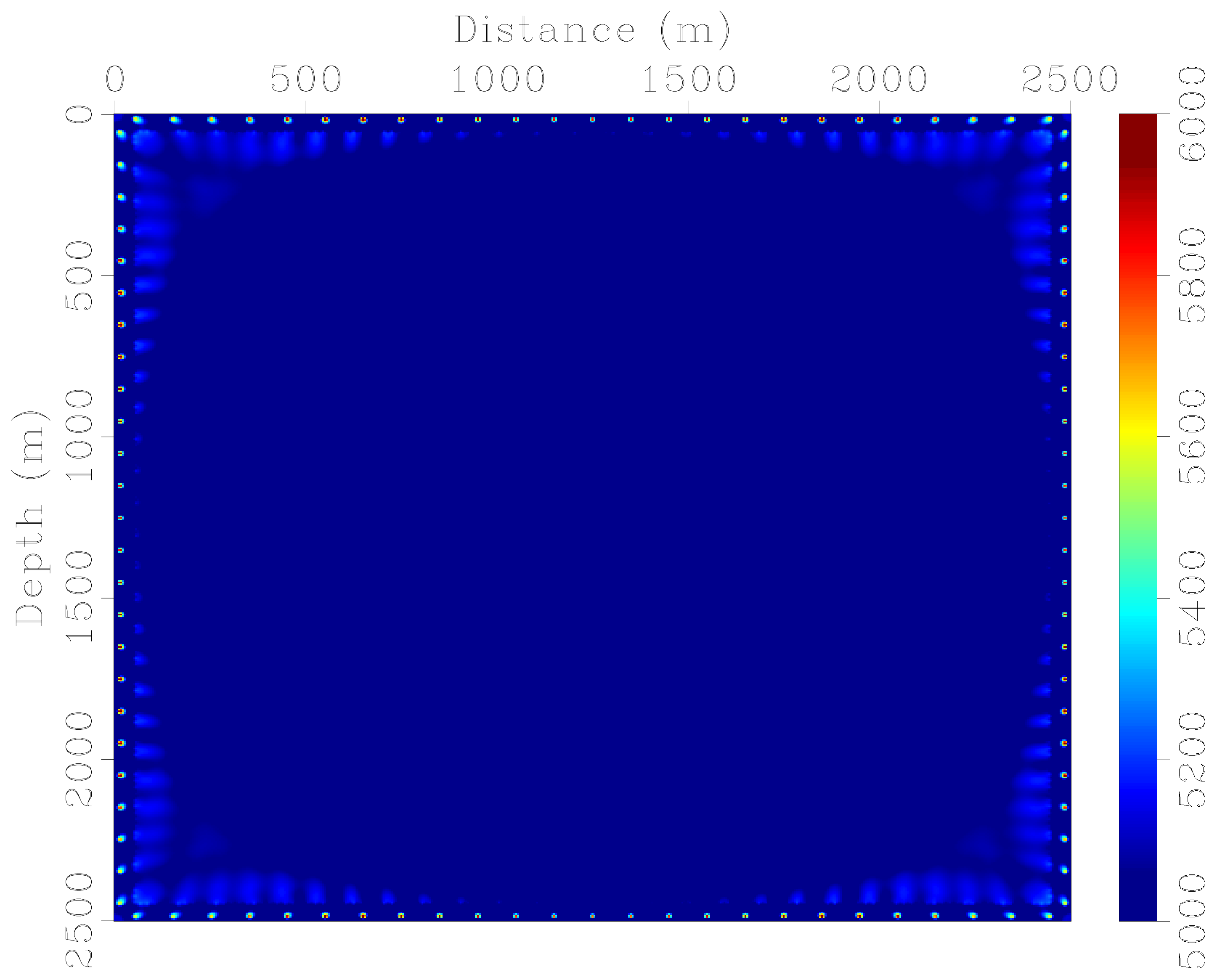} %%LPC_4s_LSQ_velsmooth81.pdf}
\end{center} 
\caption{
Inversion example 3: 
Plots of velocity models of (left) true model, 
(center) of converged model of RGLS optimization, and
(right) of converged model of LS optimization.
}
\label{fig:fully_sampled_noisy_model_convergence_behavior}
\end{figure}

Finally, we test the RGLS method with a more complicated reference velocity model shown 
in Figure \ref{fig:fully_sampled_noisy_model_convergence_behavior} (left).
The configuration of sources and receivers is the same as in inversion example 2;
sources on all sides and receivers on the 3 opposite sides.
The medium-scale details of the model are successfully recovered by RGLS optimization.
RGLS optimization stalls the data misfit after 63 iterations: the inversion then switches from RGLS to LS. (Note that the LS method is close
to the special case $\alpha = 1$ in the construction of fractionally warped data, hence the late-game switch to LS is more of a parameter
adjustment than an ad-hoc fix.)
Switching to LS is safe because observed data are now within a fraction of a wavelength of predicted data.
Using a velocity model with stronger randomness would make the RGLS method fail, because the observed data would contain many refracted
waves that could not be matched to the simple one-wave prediction by warping.

\begin{figure}
\begin{center}
\includegraphics[width=0.45\textwidth, height=0.25\textheight]{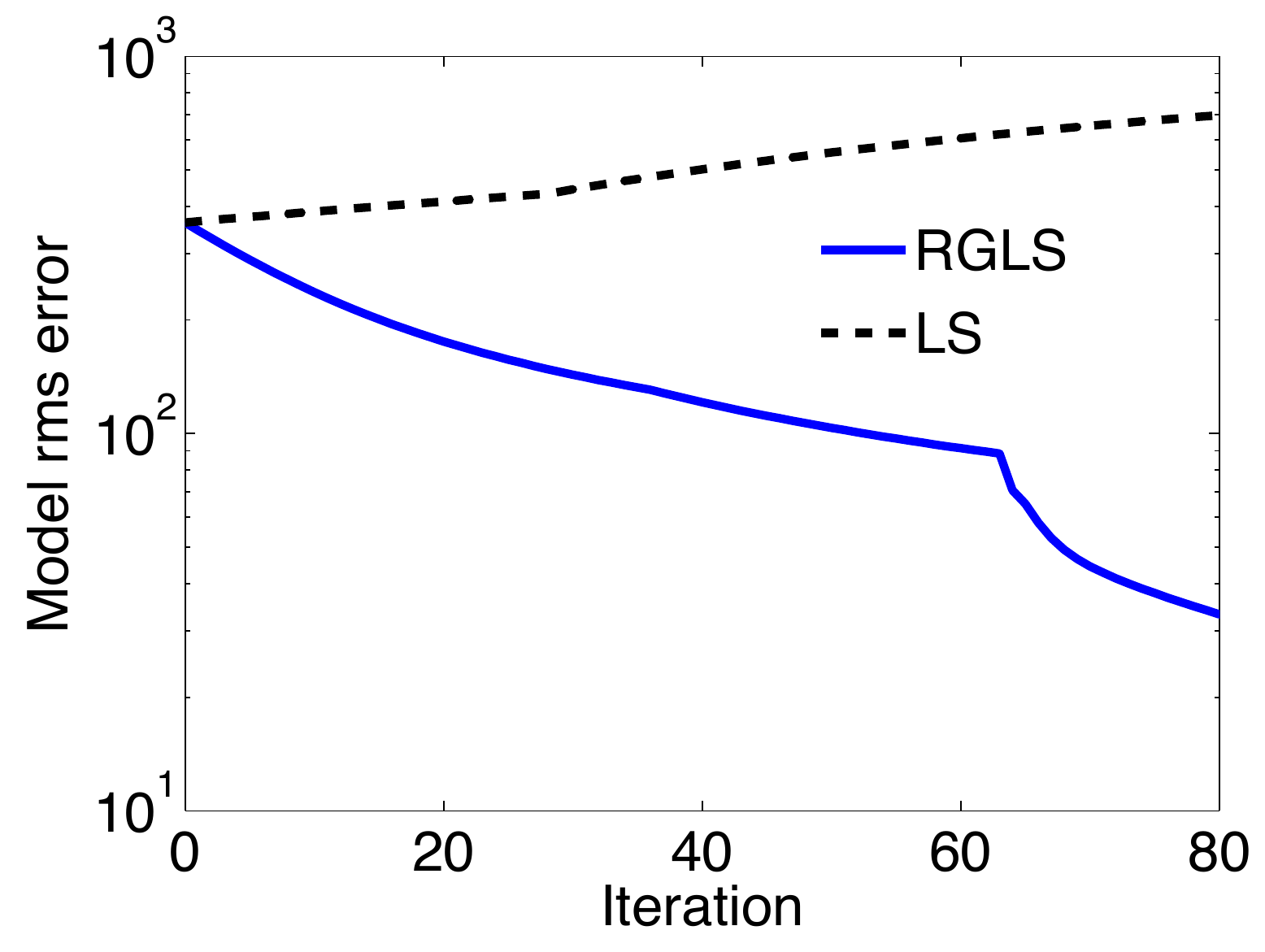} %%LPC_4s_VelRMSError.pdf}
\includegraphics[width=0.45\textwidth, height=0.25\textheight ]{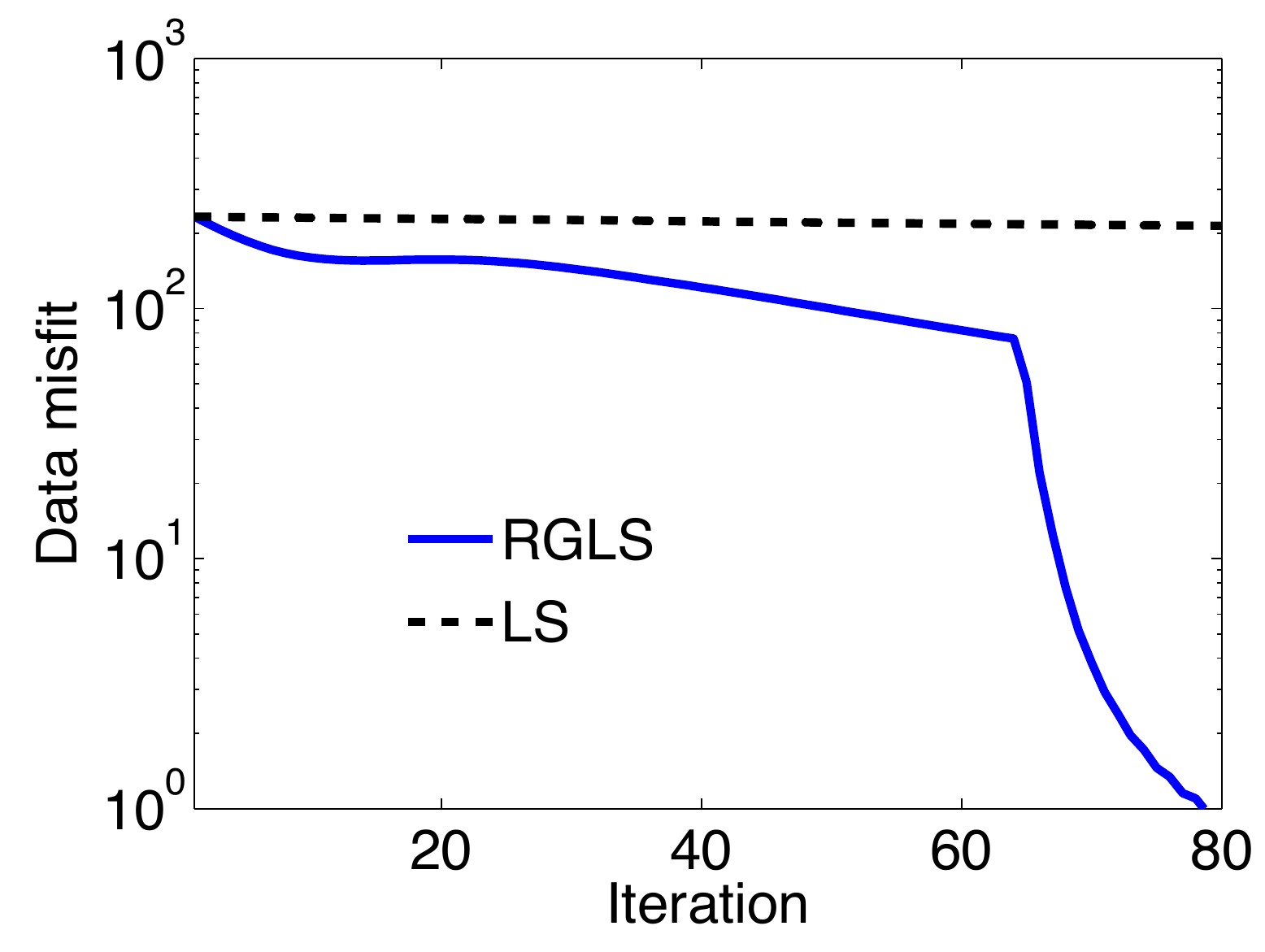} %%LPC_4s_CostDecay.pdf}
\end{center} 
\caption{Inversion example 3: 
convergence of model (velocity) rms error $V_k - V_{true} $ (left) and data misfit $J$ (right)}
\label{fig:noise_convergence_behavior}
\end{figure}

\section{Discussion}
The proposed method has a limitation; 
seismogram registration requires the predicted and observed traces to be comparable.
When the observed and predicted seismograms differ by substantially more than a warping,
the registration-based method may not make sense at all.
This issue often occurs with reflection data; it is typically difficult to create interfaces in the model that
generate the right number of reflected waves in the first iteration.
The transmission case is much more favorable as it often generates a prediction with the same number of arrivals as the data, albeit with
large traveltime delays. 

\section{Conclusions}
We present registration guided least-squares (RGLS) as a way to 
overcome the cycle-skipping problem in FWI, thereby extending the basin of attraction to the global minimizer.
The successful application of the RGLS method to seismic inversion problems is demonstrated in the transmission setting, where the
conventional LS method often converges to a wrong model.
The proposed method substitutes a transported version of the prediction, 
referred to as fractionally warped data, 
for the observed data in the conventional least-squares misfit residual.
In order to generate transported data, 
mappings in the form of piecewise polynomials 
are found through a non-convex optimization formulation.
The non-convex optimization problem is tackled in a multiscale manner 
similar to frequency sweep/continuation in frequency domain FWI.
In order to create the low frequencies which may be absent in data,
low-frequency augmented (LFA) signals are proposed and demonstrated to 
provide a satisfying alternative to the raw seismograms for the registration step.
A method using the envelope property of the Hilbert transform
is proposed for this LFA transformation.
Three inversion examples using seismogram registration and the RGLS method
show that the proposed method decreases model errors monotonically
while it allows the data misfit to increase temporarily prior to eventual convergence.

\bibliographystyle{seg}  % style file is seg.bst
\bibliography{ArxivGeoOptim}

\begin{thebibliography}{}
\itemsep0pt

\bibitem[Albertin, 2011]{Albertin2011}
Albertin, U.,  2011, An improved gradient computation for adjoint
  wave‚Äêequation reflection tomography: SEG Technical Program Expanded
  Abstracts 2011,  3969--3973.

\bibitem[Anderson and Gaby, 1983]{Anderson1983221}
Anderson, K.~R., and J.~E. Gaby,  1983, Dynamic waveform matching: Information
  Sciences, {\bf 31}, 221 -- 242.

\bibitem[Benedetto, 1997]{HarmonicAnalysis1997}
Benedetto, J.~J.,  1997, Harmonic analysis and applications: CRC-Press.

\bibitem[Bozdağ et~al., 2011]{Misfit_phase_envelope}
Bozdağ, E., J. Trampert, and J. Tromp,  2011, Misfit functions for full
  waveform inversion based on instantaneous phase and envelope measurements:
  Geophysical Journal International, {\bf 185}, 845--870.

\bibitem[Bregman et~al., 1989]{Bregman01021989}
Bregman, N.~D., R.~C. Bailey, and C.~H. Chapman,  1989, Crosshole seismic
  tomography: Geophysics, {\bf 54}, 200--215.

\bibitem[Fei and Williamson, 2010]{Fei2010}
Fei, W., and P. Williamson,  2010, On the gradient artifacts in migration
  velocity analysis based on differential semblance optimization: SEG Technical
  Program Expanded Abstracts 2010,  4071--4076.

\bibitem[Fomel and Jin, 2009]{FomelJin2009}
Fomel, S., and L. Jin,  2009, Time-lapse image registration using the local
  similarity attribute: Geophysics, {\bf 74}, A7--A11.

\bibitem[Fomel and van~der Baan, 2010]{FomelBaan2010}
Fomel, S., and M. van~der Baan,  2010, Local similarity with the envelope as a
  seismic phase detector: SEG Technical Program Expanded Abstracts 2010,
  1555--1559.

\bibitem[Haker and Tannenbaum, 2001]{938878}
Haker, S., and A. Tannenbaum,  2001, Optimal mass transport and image
  registration: Variational and Level Set Methods in Computer Vision, 2001.
  Proceedings. IEEE Workshop on, 29 --36.

\bibitem[Hale, 2013]{HaleGeoPhysics2013}
Hale, D.,  2013, Dynamic warping of seismic images: Geophysics, {\bf 78},
  S105--S115.

\bibitem[Kennett and Fichtner, 2012]{Kennett2012GJI}
Kennett, B. L.~N., and A. Fichtner,  2012, A unified concept for comparison of
  seismograms using transfer functions: Geophysical Journal International, {\bf
  191}, 1403--1416.

\bibitem[Liner and Clapp, 2004]{Liner01112004}
Liner, C.~L., and R.~G. Clapp,  2004, Nonlinear pairwise alignment of seismic
  traces: The Leading Edge, {\bf 23}, 1146--1150.

\bibitem[Luo and Schuster, 1991]{luo645}
Luo, Y., and G.~T. Schuster,  1991, Wave-equation traveltime inversion:
  Geophysics, {\bf 56}, 645--653.

\bibitem[Maggi et~al., 2009]{Maggi2009}
Maggi, A., C. Tape, M. Chen, D. Chao, and J. Tromp,  2009, An automated
  time-window selection algorithm for seismic tomography: Geophysical Journal
  International, {\bf 178}, 257--281.

\bibitem[Plessix, 2009]{Plessix2009}
Plessix, R.-E.,  2009, Three-dimensional frequency-domain full-waveform
  inversion with an iterative solver: Geophysics, {\bf 74}, WCC149--WCC157.

\bibitem[Pratt, 1999]{Pratt99seismicwaveform}
Pratt, R.~G.,  1999, Seismic waveform inversion in the frequency domain, part
  1: Theory and verification in a physical scale model: Geophysics, {\bf 64},
  888--901.

\bibitem[Pratt and Goulty, 1991]{Pratt01021991}
Pratt, R.~G., and N.~R. Goulty,  1991, Combining wave-equation imaging with
  traveltime tomography to form high-resolution images from crosshole data:
  Geophysics, {\bf 56}, 208--224.

\bibitem[Prieux et~al., 2012]{GPR1099}
Prieux, V., G. Lambar\'e, S. Operto, and J. Virieux,  2012, Building starting
  models for full waveform inversion from wide-aperture data by
  stereotomography: Geophysical Prospecting,  no--no.

\bibitem[Ramsay and Li, 1998]{ramsay1998}
Ramsay, J.~O., and X. Li,  1998, Curve registration: Journal of the Royal
  Statistical Society. Series B, {\bf 60}, 351--363.

\bibitem[Sakoe and Chiba, 1978]{saoke78}
Sakoe, H., and S. Chiba,  1978, Dynamic programming algorithm optimization for
  spoken word recognition: IEEE Transactions on Acoustics, Speech, and Signal
  Processing,  43--49.

\bibitem[Sava and Biondi, 2004]{sava_biondi2004geoprospect}
Sava, P., and B. Biondi,  2004, Wave-equation migration velocity analysis. {I}.
  theory: Geophysical Prospecting, {\bf 52}, 593--606.

\bibitem[Tarantola and Valette, 1982]{Tarantola1982}
Tarantola, A., and B. Valette,  1982, Generalized nonlinear inverse problems
  solved using the least squares criterion: Rev. Geophys., {\bf 20}, 219--232.

\bibitem[Tromp et~al., 2005]{Tromp2005}
Tromp, J., C. Tape, and Q. Liu,  2005, Seismic tomography, adjoint methods,
  time reversal and banana-doughnut kernels: Geophysical Journal International,
  {\bf 160}, 195--216.

\bibitem[Virieux and Operto, 2009]{Virieux2009}
Virieux, J., and S. Operto,  November-December 2009, An overview of
  full-waveform inversion in exploration geophysics: Geophysics, {\bf 74},
  WCC1--WCC26.

\end{thebibliography}

\end{document}